\documentclass{article}
\usepackage{epsfig,ulem,amssymb}
\usepackage[margin=0.94in]{geometry}
\usepackage{hyperref}
\usepackage{graphicx}
\usepackage{epstopdf}
\usepackage{amsmath}
\usepackage{bm}
\usepackage{setspace}
\usepackage{mathrsfs}
\usepackage{subcaption,comment}
\usepackage[dvipsnames]{xcolor}
\usepackage{mathtools}
\usepackage{multirow}
\usepackage{lipsum}
\usepackage{enumerate}
\usepackage[affil-it]{authblk}
\usepackage{cite}
\usepackage{leftidx}
\usepackage{relsize}
\usepackage{braket}
\usepackage{relsize}
\usepackage{float}
\DeclarePairedDelimiter{\ceil}{\lceil}{\rceil}
\usepackage[autostyle]{csquotes}

\usepackage{fontenc}
\usepackage[right]{lineno}

\providecommand{\keywords}[1]{\textbf{\textit{Keywords---}} #1}

\title{A Ritz-based Finite Element Method for a Fractional-Order Boundary Value Problem of Nonlocal Elasticity}
\author[1]{Sansit Patnaik$^\dagger$}
\affil[1]{School of Mechanical Engineering, Ray W. Herrick Laboratories, Purdue University, West Lafayette, IN 47907}
\author[1]{Sai Sidhardh}
\author[1]{Fabio Semperlotti$^\dagger$}

\begin{document}
\date{}
\maketitle

\begin{abstract}
We present the analytical formulation and the finite element solution of a fractional-order nonlocal continuum model of a Euler-Bernoulli beam. Employing consistent definitions for the fractional-order kinematic relations, the governing equations and the associated boundary conditions are derived based on variational principles. Remarkably, the fractional-order nonlocal model gives rise to a self-adjoint and positive-definite system accepting a unique solution. Further, owing to the difficulty in obtaining analytical solutions to this boundary value problem, a finite element model for the fractional-order governing equations is presented. Following a thorough validation with benchmark problems, the fractional finite element model (f-FEM) is used to study the nonlocal response of a Euler-Bernoulli beam subjected to various loading and boundary conditions. The fractional-order positive definite system will be used here to address some paradoxical results obtained for nonlocal beams through classical integral approaches to nonlocal elasticity. Although presented in the context of a 1D Euler-Bernoulli beam, the f-FEM formulation is very general and could be extended to the solution of any general fractional-order boundary value problem.\\

\noindent\keywords{Fractional Calculus, Nonlocal Beams, Variational Calculus, Finite Element Method}\\
\noindent$^\dagger$ All correspondence should be addressed to: \textit{spatnai@purdue.edu} or \textit{fsemperl@purdue.edu}

\end{abstract}

\section{Introduction}
\label{sec:Introduction}
Recent theoretical and experimental studies have shown that scale-dependent effects are prominent in the response of several structures ranging from layered and porous media \cite{gurevich1995velocity,biot1956theory,buonocore2018occurrence}, to random and fractal media \cite{szabo1994time,fellah2004verification}, to media with damage and cracks \cite{stulov2016frequency,berkowitz1995characterization}, to biomedical materials like tissues and bones \cite{szabo1994time,adbone}. Size-dependent effects are particularly prominent in microstructures and nanostructures that have far-reaching applications in atomic devices, micro/nano-electromechanical devices and sensors \cite{sumelka2015fractional,civalek2016simple,rahimi2017linear}. These devices are primarily made from a combination of beams or other slender structures like plates and shells. Accurate modeling of the response of these structures is paramount in many engineering as well as biomedical applications.

The coexistence of different spatial scales in the above mentioned classes of structural problems renders the response fully nonlocal \cite{kroner1967elasticity,eringen1972linear,eringen1972nonlocal}. The inability of the classical continuum theory in capturing scale effects prevents its use in these types of applications and fostered the interest in the development of nonlocal continuum theories. 
Nonlocal continuum theories enrich the classical (local) governing equations describing the response at a point with information of the behaviour of points contained in a prescribed area of influence. The key principle behind these nonlocal continuum theories is that, all the particles located inside this area, also known as the horizon of influence or horizon of nonlocality, influence one another by means of long range cohesive forces \cite{kroner1967elasticity,eringen1972linear,eringen1972nonlocal}. Seminal works from Kro\"ner \cite{kroner1967elasticity}, Eringen \cite{eringen1972nonlocal}, and several other authors \cite{nowinski1984nonlocal,nowinski1986non} have explored the role of nonlocality in elasticity and laid its theoretical foundation. Further, several theories based on integral methods, gradient methods and very recently, the peridynamic approach have been developed to capture these long range energy exchanges and analyze their effect on the response of structures.
Gradient elasticity theories \cite{peerlings2001critical,aifantis2003update,guha2015review} account for the nonlocal behavior by introducing strain gradient dependent terms in the stress-strain constitutive law. These strain gradient theories have been extensively used in the study of different structures \cite{sidhardh2019exact,sidhardh2018element,sidhardh2018inclusion,guruprasad2019some}. Integral methods \cite{polizzotto2001nonlocal,bavzant2002nonlocal,sidhardh2018effect} model nonlocal effects by defining the constitutive law in the form of a convolution integral between the strain and the spatially dependent elastic properties over the horizon of nonlocality. Recently, Silling \cite{silling2000reformulation} proposed the peridynamic approach as an alternative theory that is better suited to model structures involving a dynamic evolution of discontinuities.

In recent years, fractional calculus has emerged as a powerful mathematical tool to model a variety of nonlocal and multiscale phenomena. Fractional derivatives, which are a differ-integral class of operators, are intrinsically multiscale and provide a natural way to account for nonlocal effects. As a result, time-fractional operators enable memory effects (i.e. the response of a system is a function of its past history) while space-fractional operators can account for nonlocal and scale effects. These characteristics of fractional operators have led to a surge of interest in fractional operators and in their applications to the simulation of several physical problems. Areas that have seen the largest number of applications include the formulation of constitutive equations for viscoelastic materials \cite{bagley1983theoretical,koeller1984applications,chatterjee2005statistical}, model-order reduction of lumped parameter systems \cite{hollkamp2018model}, and modeling of transport processes in complex media \cite{mainardi1996fractional,benson2000application,hollkamp2019analysis,hollkamp2020application}.

Given the multiscale nature of fractional operators, fractional calculus has also found wide-spread application in nonlocal elasticity. 
Riesz-type fractional derivatives have been shown to emerge as the continuum limit of discrete systems (e.g. such as chains and lattices) with power-law long-range interactions \cite{tarasov2013lattice,tarasov2013review}. 
Space-fractional derivatives have been used to formulate nonlocal constitutive laws \cite{drapaca2012fractional,carpinteri2014nonlocal,sumelka2014thermoelasticity,Sumelka1,sumelka2016fractional} as well as to account for microscopic interaction forces \cite{lazopoulos2006non,cottone2009fractional,di2008long}. Space-fractional derivatives have been employed to capture attenuation including a variety of conditions such as interatomic nonlocal forces \cite{di2008long,tarasov2013lattice,tarasov2013review}, nonlocal stress-strain constitutive relations \cite{atanackovic2009generalized}, and even bandgaps in periodic media \cite{hollkamp2019analysis}. Very recently fractional-order nonlocal theories have been extended to model and analyze the static response, buckling characteristics, as well as the dynamic response of nonlocal beams \cite{sumelka2015fractional,rahimi2017linear,alotta2017finite}.

In this study, we build upon the fractional-order nonlocal continuum model proposed in \cite{sumelka2014thermoelasticity,sumelka2016fractional,sumelka2016geometrical} to develop a fractional-order constitutive relation for the Euler-Bernoulli beam. The fractional-order nonlocal continuum model proposed in \cite{sumelka2014thermoelasticity,sumelka2016fractional,sumelka2016geometrical} is shown to be frame-invariant, dimensionally consistent, and, unlike other fractional-order nonlocal theories, requires integer-order boundary conditions that accept a clear physical interpretation. However, we emphasize that the definition of the fractional deformation tensor adopted in this study is different from that used in \cite{sumelka2014thermoelasticity,sumelka2016fractional,sumelka2016geometrical} which has important implications on the resulting fractional-order framework. The overall goal of this study is two-fold. First, we derive the governing equations for the nonlocal beam in a strong form using variational principles. More specifically, the governing equations are derived by minimization of the total potential energy of the beam. This approach is different from that proposed in \cite{sumelka2015fractional,rahimi2017linear,alotta2017finite}, where the equations of motion describing the nonlocal beam had been derived using Newton's approach of force and moment equilibrium and from different fractional-order nonlocal constitutive relationships.
Additionally, we show that the fractional-order modeling of the nonlocal beam results in a self-adjoint system with a quadratic potential energy, irrespective of the boundary conditions. This result is in sharp contrast with the integral nonlocal methods available in the literature for which it is not possible to define a self-adjoint quadratic potential energy \cite{reddy2007nonlocal,challamel2008small,challamel2014nonconservativeness}.
Second, we formulate a fully consistent and highly accurate fractional-order finite element method (f-FEM) to numerically investigate the response of the fractional-order nonlocal beam. Although several FE formulations for fractional-order equations have been proposed in the literature, they are based on Galerkin or Petrov-Galerkin methods that are capable of solving hyperbolic and parabolic differential equations involving transport processes \cite{agrawal2008general,deng2008finite,zheng2010note,zhang2010galerkin,jiang2011high,bu2014galerkin,liu2015two,jin2015galerkin,wang2015petrov,yang2017finite}. We develop a Ritz FEM that is capable of obtaining the numerical solution of the fractional-order elliptic boundary value problem (BVP) that describes the static response of the fractional-order nonlocal beam. Although Ritz FEMs for classical nonlocal elasticity problems have been developed, they do not extend to fractional-order nonlocal modeling, because the attenuation function capturing the nonlocal interactions in the fractional-order model involves a singularity within the kernel \cite{podlubny1998fractional}. In this study, we have outlined a strategy to treat this singularity in the fractional derivatives. Further, by using the f-FEM we show that, independently from the boundary conditions, the fractional-order theory predicts a softening behaviour for the fractional-order beam as the nonlocality degree increases. With these results, we explain the paradoxical predictions of hardening and absence of nonlocal effects for certain combinations of boundary conditions, as predicted by classical integral approaches to nonlocal elasticity \cite{challamel2008small,challamel2014nonconservativeness,khodabakhshi2015unified}.

The remainder of the paper is structured as follows: first, we present the fractional-order model of a nonlocal Euler-Bernoulli beam. Next, we derive the governing equations of the beam in strong form using variational principles. Further we derive a strategy for obtaining the numerical solution to the beam governing equation using fractional-order FEM. Finally we validate the fractional-order FEM, establish its convergence, and then use it to analyze the effect of the fractional-order nonlocality on the static response of the beam under different types of loading conditions.

\section{Nonlocal Euler-Bernoulli Beam Model}
\label{sec:NBM}
Previous works conducted on the development of nonlocal continuum theories based on fractional calculus have highlighted its ability to combine the strengths of both gradient and integral based methods while at the same time addressing a few important shortcomings of these integer order formulations \cite{drapaca2012fractional,carpinteri2014nonlocal,sumelka2014thermoelasticity}. Gradient elasticity theories provide a satisfactory description of the material micro structure, but they introduce serious difficulties when enforcing the boundary conditions associated with the strain gradient-dependent terms \cite{peerlings2001critical,aifantis2003update}. 
On the other side, the integral methods are better suited to deal with boundary conditions but require the attenuation functions to have a positive Fourier transform everywhere in order to avoid instabilities \cite{bavzant1984instability,bavzant2002nonlocal}. To this regard, note that the kernel used in fractional derivatives is positive everywhere \cite{podlubny1998fractional}. Thus, by formulating the constitutive relations using space-fractional derivatives, the resulting nonlocal theory effectively combines features characteristic of both gradient-based and integral-based methods. Unlike gradient elasticity methods, additional boundary conditions are not required when using Caputo fractional derivatives \cite{hollkamp2019analysis}. The nonlocal beam theory presented in this work builds on the fractional-order nonlocal continuum formulation presented in \cite{sumelka2014thermoelasticity,sumelka2015fractional,sumelka2016fractional} where nonlocality is accounted for by means of a fractional-order deformation gradient tensor. However, as highlighted in the introduction, the definition of the fractional deformation tensor adopted in this study is different from that used in \cite{sumelka2014thermoelasticity,sumelka2015fractional,sumelka2016fractional}. In the following we present and discuss important aspects of the fractional-order continuum theory and then use the same to develop the fractional-order Euler-Bernoulli beam theory.

\subsection{Nonlocal Continuum Formulation}
\label{sec:nonlocal_continuum_formulation}
\begin{figure}[h]
	\centering
	\includegraphics[width=\linewidth]{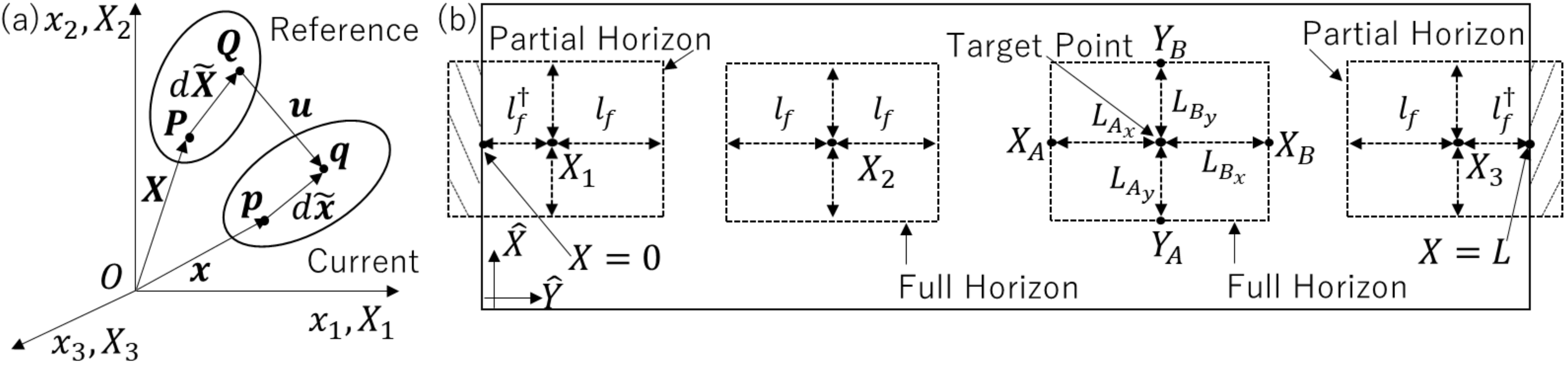}
	\caption{\label{fig1} (a) Schematic indicating the infinitesimal material $\mathrm{d}\tilde{\textbf{X}}$ and spatial $\mathrm{d}\tilde{\textbf{x}}$ line elements in the nonlocal medium subject to the displacement field $\textbf{u}$. (b) Horizon of nonlocality and length scales at three different material points $\textbf{X}_1$, $\textbf{X}_2$, and $\textbf{X}_3$ in a 2D domain. Note from Fig.~(\ref{fig1}) that in the $\hat{X}$ direction, $\textbf{X}_2$ has a horizon of nonlocality equal to $l_f$ on both the left and the right sides, while the horizon of nonlocality at the points $\textbf{X}_1$ and $\textbf{X}_3$ are truncated to $l_f^\dagger$ such that $l_f^\dagger<l_f$, on the left and the right sides, respectively. Clearly, a similar approach is applied also to account for boundaries in the $\hat{Y}$ direction. The nonlocal model can account for a partial (i.e. asymmetric) horizon condition that occurs at points ${\textbf{X}}$ close to a boundary or interface.}
\end{figure}

We perform the deformation analysis of a nonlocal solid by introducing two stationary configurations, namely, the reference (undeformed) and the current (deformed), in analogy with the traditional continuum approach to mechanics. A motion of the body from the reference configuration (denoted as $\textbf{X}$) to the current configuration (denoted as $\textbf{x}$) is considered:
\begin{equation}
\label{motion_description}
\textbf{x}=\bm{\Psi}(\textbf{X})
\end{equation}
such that $\bm{\Psi}(\textbf{X})$ is a continuous and invertible mapping operation.
The relative position of two point particles located at $P$ and $Q$ in the reference configuration of the nonlocal medium is denoted by ${\mathrm{d}\tilde{\textbf{X}}}$ (see Fig.~(\ref{fig1})). After deformation due to the motion $\bm{\Psi}(\textbf{X})$, the particles move to new positions $p$ and $q$, such that the relative position vector between them is ${\mathrm{d}\tilde{\textbf{x}}}$. Thus $\mathrm{d}\tilde{\textbf{X}}$ and $\mathrm{d}\tilde{\textbf{x}}$ are the material and spatial differential line elements in the nonlocal medium (conceptually analogous to the classical differential line elements $\mathrm{d}{\textbf{X}}$ and $\mathrm{d}{\textbf{x}}$).

In the classical continuum formulation, the differential line elements in the reference and current configurations (that is $\mathrm{d}{\textbf{X}}$ and $\mathrm{d}{\textbf{x}}$) are related using the classical integer-order deformation gradient tensor as:
\begin{equation}
\label{classical_F}
\mathrm{d}\textbf{x}=\big[D^{1}_\textbf{X}\bm{\Psi}(\textbf{X})\big]\mathrm{d}\textbf{X}=[\textbf{F}(\textbf{X})]\mathrm{d}\textbf{X}
\end{equation}
where $D^{1}_\textbf{X}(\cdot)$ denotes the first integer-order spatial derivative with respect to the reference coordinates. In the fractional-order formulation, the differential line elements of the nonlocal medium are modeled by imposing a fractional-order transformation on the classical differential line elements as follows:
\begin{subequations}
\label{Fractional_F}
\begin{equation}
\mathrm{d}\tilde{\textbf{x}}=\big[D^{\alpha}_\textbf{X}\bm{\Psi}(\textbf{X})\big]\mathrm{d} \textbf{X}=\big[\tilde{\textbf{F}}_{X}(\textbf{X})\big]\mathrm{d}\textbf{X}
\end{equation}
\begin{equation}
\mathrm{d}\tilde{\textbf{X}}=\big[D^{\alpha}_\textbf{x}\bm{\Psi}^{-1}(\textbf{x})\big]\mathrm{d} \textbf{x}=\big[\tilde{\textbf{F}}_{x}(\textbf{x})\big]\mathrm{d}\textbf{x}
\end{equation}	
\end{subequations}
where $D^{\alpha}_\square\square$ is a space-fractional derivative whose details will be presented below. Given the differ-integral nature of the space-fractional derivative, the differential line elements $\mathrm{d}\tilde{\textbf{X}}$ and $\mathrm{d}\tilde{\textbf{x}}$ are nonlocal in nature.
Using the definitions for $\mathrm{d}\tilde{\textbf{X}}$ and $\mathrm{d}\tilde{\textbf{x}}$, the fractional deformation gradient tensor $\overset{\alpha}{\textbf{F}}$ with respect to the nonlocal coordinates is obtained as:
\begin{equation}
\label{Fractional_F_net}
\frac{\mathrm{d}\tilde{\textbf{x}}}{\mathrm{d}\tilde{\textbf{X}}}=\mathop{\textbf{F}}^{\alpha}=\tilde{\textbf{F}}_{X}\textbf{F}^{-1}\tilde{\textbf{F}}_{x}^{-1}
\end{equation}

The space-fractional derivative $D^{\alpha}_\textbf{X}\bm{\Psi}(\textbf{X},t)$ is taken according to a Riesz-Caputo (RC) definition with order $\alpha\in(0,1)$ defined on the interval $\textbf{X} \in (\textbf{X}_A,\textbf{X}_B) \subseteq \mathbb{R}^3 $ and is given by:
\begin{subequations}
	\label{RC_definition}
	\begin{equation}
	D^{\alpha}_\textbf{X}\bm{\Psi}(\textbf{X},t)=\frac{1}{2}\Gamma(2-\alpha)\big[\textbf{L}_{A}^{\alpha-1}~ {}^C_{\textbf{X}_{A}}D^{\alpha}_{\textbf{X}} \bm{\Psi}(\textbf{X},t) - \textbf{L}_{B}^{\alpha-1}~ {}^C_{\textbf{X}}D^{\alpha}_{\textbf{X}_{B}}\bm{\Psi}(\textbf{X},t)\big]
	\end{equation}
	\begin{equation}
	D^{\alpha}_{X_j}\Psi_i(\textbf{X},t)=\frac{1}{2}\Gamma(2-\alpha)\big[L_{A_j}^{\alpha-1}~{}^C_{X_{A_j}}D^{\alpha}_{X_j} \Psi_i(\textbf{X},t) - L_{B_j}^{\alpha-1}~{}^C_{X_j}D^{\alpha}_{X_{B_j}} {\Psi_i}(\textbf{X},t)\big]
	\end{equation}
\end{subequations}
where, ${}^C_{\textbf{X}_{A}}D^{\alpha}_{\textbf{X}}\bm{\Psi}(\textbf{X},t)$ and ${}^C_{\textbf{X}}D^{\alpha}_{\textbf{X}_{B}}\bm{\Psi}(\textbf{X},t)$ are the left- and right-handed Caputo derivatives of $\bm{\Psi}(\textbf{X},t)$ respectively, and $\Gamma(\cdot)$ is the Gamma function. In the indicial expression in Eq.~(\ref{RC_definition}b), $L_{A_j}$ and $L_{B_j}$ are length scales along the $j^{th}$ direction in the reference configuration. 
The index $j$ in Eq.~(\ref{RC_definition}b) is not a repeated index because the length scales are scalar multipliers. In the current configuration, these length scales are denoted as $l_{A_j}$ and $l_{B_j}$. The interval of the fractional derivative $(\textbf{X}_A,\textbf{X}_B)$ defines the horizon of nonlocality (also called attenuation range in classical nonlocal elasticity) which is schematically shown in Fig.~(\ref{fig1}) for a generic point $\textbf{X}\in\mathbb{R}^2$.
In other terms, it defines the set of all points in the solid that influence the elastic response at $\textbf{X}$ (or, equivalently, the characteristic distance beyond which information of nonlocal interactions is no longer accounted for in the derivative). We have shown in Appendix 1 that, for a frame-invariant model, it is required that the length scales $\textbf{L}_A=\textbf{X}-\textbf{X}_A$ and $\textbf{L}_B=\textbf{X}_B-\textbf{X}$. Hence, it follows that the length scales, $L_{A_j}$ and $L_{B_j}$, physically denote the dimension of the horizon of nonlocality to the left and right of point $\textbf{X}$ along the $j^{th}$ direction. The length scales have been schematically illustrated in Fig.~(\ref{fig1}b).
Further, the length scale parameters $L_{A_j}^{\alpha-1}$ and $L_{B_j}^{\alpha-1}$ ensure the dimensional consistency of the deformation gradient tensor.

We emphasize that while \cite{sumelka2016fractional} focused primarily on anisotropic nonlocality, the model presented in this study is for isotropic materials. The introduction of the different length scales ($\textbf{L}_{A}$ and $\textbf{L}_B$) is to enable the formulation to deal with possible asymmetries in the horizon of nonlocality (e.g. resulting from a truncation of the horizon when approaching a boundary or an interface).
More specifically, the different length scales enable an efficient and accurate treatment of the frame invariance and ensure a completeness of the kernel (established in the following) in the presence of asymmetric horizons, material boundaries, and interfaces (Fig.~(\ref{fig1}b)).

The definition of the RC fractional derivative in Eq.~(\ref{RC_definition}) ensures the completeness of the power-law convolution kernel within the fractional-order derivative. Note that the lower terminal is $\textbf{X}_A=\textbf{X}-\textbf{L}_A$ and the upper terminal is $\textbf{X}_B=\textbf{X}+\textbf{L}_B$. This definition allows the length scales $\textbf{L}_A$ and $\textbf{L}_B$ to be truncated when the point $\textbf{X}$ approaches a boundary (see Fig.~(\ref{fig1}b)). In other terms, the fractional-order model defined in this study
allows for position dependent length scales $\textbf{L}_A(\textbf{X})$ and $\textbf{L}_B(\textbf{X})$ for points where the horizon is truncated due to the presence of boundaries or any other source of discontinuities.
It follows that the terminals of the RC derivative are properly modified hence resulting in a complete kernel over the truncated domain.

The completeness of the kernel can also be illustrated by investigating the nature of the fractional-order model at points on the boundary. For a 3D fractional-order continuum model, we investigate the fractional deformation gradient tensor for points located on the boundary, that is when either $L_{A_j}\rightarrow0$ or $L_{B_j}\rightarrow0$.
For a material point (say $\textbf{X}_0$) located on one of the boundaries (identified by the normal in the $j^{th}$ direction), one of the length scales goes to zero. Thus, the interval length of either the left-handed or the right-handed Caputo derivative is $0$. In the special case of a corner point, both lengths will go to zero simultaneously. This singularity in the RC derivative is treated by taking the interval of the left-handed (or right-handed) Caputo derivative as $\varepsilon=X_{0_j}-X_{A_j}=L_{A_j}$ (or $\varepsilon=X_{B_j}-X_{0_j}=L_{B_j}$) and evaluating the limiting condition $\varepsilon\rightarrow0$. Here below we present the expression for the RC derivative when $L_{A_j}=0$. Similar expressions hold when $L_{B_j}=0$ and for the deformed configuration ($l_{A_j}=0$ or $l_{B_j}=0$). The limiting case gives:
\begin{equation}
\label{sq10}
\lim_{L_{A_j}\to0}D^{\alpha}_{X_j}\psi_i(\textbf{X},t) = \lim_{\varepsilon\to0} \frac{(1 -\alpha)}{2} \bigg[
{\varepsilon^{\alpha-1}} \int_{X_{0_j}-\varepsilon}^{X_{0_j}} \frac{D^1_{S_j} \psi_{i}(\textbf{S},t)}{(X_j-S_j)^{\alpha}} \mathrm{d}S_j + 
{L_{B_{j}}^{\alpha-1}} \int_{X_{0_j}}^{X_{B_j}} \frac{D^1_{S_j} \psi_{i}(\textbf{S},t)}{(S_j-X_j)^{\alpha}} \mathrm{d}S_j \bigg]
\end{equation}
where $S_j$ is a dummy vector variable used to carry out the spatial convolution integral within the definition of the fractional derivative. Note that in the left-handed Caputo derivative $X_{0_j}-\varepsilon<S_j<X_{0_j}$. Since the interval length of the left-handed derivative ($=\varepsilon$) is very small, $D^1_{S_j} \psi_{i}(\textbf{S},t)$ can be assumed to be constant and equal to the boundary condition. Thus, for the left-handed Caputo derivative:
\begin{equation}
\label{sq11}
D^1_{S_j} \psi_{i}(\textbf{S},t)\big|_{\textbf{S}_0}=\frac{\mathrm{d}\psi_i(\textbf{S},t)}{\mathrm{d}S_j}\Big|_{\textbf{S}_0}=\frac{\mathrm{d}\psi_i(\textbf{X},t)}{\mathrm{d}Xj}\Big|_{\textbf{X}_0}
\end{equation}
Substituting Eq.~(\ref{sq11}) in Eq.~(\ref{sq10}) leads to the following:
\begin{equation}
\label{sq12}
\begin{split}
\lim_{L_{A_j}\to0} D^{\alpha}_{X_j}\psi_i(\textbf{X},t) = \lim_{\varepsilon\to0}\frac{(1-\alpha)}{2}
\biggl[ {\varepsilon^{\alpha-1}} \frac{\mathrm{d}\psi_i(\textbf{X},t)}{\mathrm{d}Xj}\Big|_{\textbf{X}_0} \int_{X_{0_j}-\varepsilon}^{X_{0_j}}\frac{\mathrm{d}S_j}{(X_j-S_j)^{\alpha}} 
+ {L_{B_{j}}^{\alpha-1}} \int_{X_{0_j}}^{X_{B_j}} \frac{D^1_{S_j} \psi_{i}(\textbf{S},t)}{(S_j-X_j)^{\alpha}} \mathrm{d}S_j \biggr]
\end{split}
\end{equation}
Note that the first integral in the above equation is equal to $\varepsilon^{1-\alpha}$ which cancels out the singular term $\varepsilon^{\alpha-1}$. It follows that:
\begin{equation}
\label{sq13}
\lim_{L_{A_j}\to 0} D^{\alpha}_{X_j} \psi_i (\textbf{X},t)= \frac{1}{2} \left[ \frac{\mathrm{d}\psi_i(\textbf{X},t)}{\mathrm{d}X_j} \bigg|_{\textbf{X}_0} + (1-\alpha) L_{B_{j}}^{\alpha-1} \int_{X_{0_j}}^{X_{B_j}} \frac{D^1_{S_j}\psi_{i}(\textbf{S},t)}{(S_j-X_j)^{\alpha}} dS_j \right]
\end{equation}
From Eq.~(\ref{sq13}) it is immediate to observe that while the right-handed Caputo derivative captures nonlocality ahead of the point $X_0$ (in the $j^{th}$ direction), the left-handed derivative is reduced to the classical first-order derivative. This suggests that the truncation of the nonlocal horizon (and the corresponding convolution) at the boundary has been accounted for in a consistent manner.

Note that the fractional-order formulation presented above also admits the fractional-order $\alpha\in(0,1)$ as a parameter. The order $\alpha$ characterizes the strength of the nonlocal interaction over the spatial interval $(\textbf{X}_A,\textbf{X}_B)$. The power law kernel $1/|\textbf{X}|^{\alpha}$ embedded in the definition of the fractional derivatives is analogous to the attenuation function commonly used in classical integral theories of nonlocal elasticity. For $\alpha$ close to $1.0$, the power-law kernel behaves analogous to a dirac-delta function, and reduces the model to be purely local. However, for values of $\alpha$ increasingly smaller than $1.0$, the contribution of points distant from the target point $\textbf{X}$ plays a significant role in the response at $\textbf{X}$, thereby accounting for the effect of nonlocal (long-range) interactions into the model.

We also briefly discuss the specific impact of the order and the length-scales on the degree of nonlocality. Note that the convolution kernel used within the definition of the fractional-order derivative is a power-law kernel which is monotonically decreasing in nature. It appears that, by reducing the value of the order $\alpha$, the magnitude of this kernel increases. More specifically, for a fixed point $x$ interacting with a specific point $s$ in its horizon of nonlocality, $\left(\frac{1}{|\textbf{X}-\textbf{S}|}\right)^{\alpha_1} > \left(\frac{1}{|\textbf{X}-\textbf{S}|}\right)^{\alpha_2} \forall~\alpha_1 < \alpha_2$. Since, the magnitude of the kernel increases, the strength of the corresponding nonlocal interactions increases and consequently, the degree of nonlocality increases.
Similarly, by increasing the value of the length scales, the size of the horizon of nonlocality increases. Hence, information corresponding to the nonlocal interactions between a larger number of points within the solid are accounted for in the formulation. Thus the degree of nonlocality increases.

In analogy with the classical strain measures, the nonlocal strain is defined using the nonlocal fractional-order differential line elements as $\mathrm{d}\tilde{\textbf{x}} \mathrm{d}\tilde{\textbf{x}} - \mathrm{d}\tilde{\textbf{X}}\mathrm{d}\tilde{\textbf{X}}$. The fractional deformation gradient tensor has been used to derive the infinitesimal fractional-order nonlocal strain tensor as:
\begin{equation}
\label{infinitesimal_fractional_strain}
\tilde{\bm{\epsilon}}=\frac{1}{2}\bigr(\nabla^\alpha {\textbf{U}}_X+\nabla^\alpha {\textbf{U}}_X^{T}\bigl)=\frac{1}{2}\bigr(\nabla^\alpha {\textbf{u}}_x+\nabla^\alpha {\textbf{u}}_x^{T}\bigl)
\end{equation}
where $\textbf{U}(\textbf{X})=\textbf{x}(\textbf{X})-\textbf{X}$ and $\textbf{u}(\textbf{x})=\textbf{x}-\textbf{X}(\textbf{x})$ are the displacement fields in the Lagrangian and Eulerian coordinates, respectively. In Eq.~(\ref{infinitesimal_fractional_strain}), $\nabla^\alpha\textbf{U}_X$ is the fractional gradient given as $\nabla^\alpha\left[\textbf{U}_{X}\right]_{ij}=D^{\alpha}_{X_j}U_i$.

Further, stress in the nonlocal isotropic medium is given analogous to the local isotropic case as:
\begin{equation}
\label{stress_equation}
\tilde{\sigma}_{ij}=2\mu\tilde{\varepsilon}_{ij}+\lambda\tilde{\bm{\varepsilon}}_{kk}\delta_{ij}
\end{equation}
where $\lambda$ and $\mu$ are Lam\'{e} constants. As expected, classical continuum mechanics relations are recovered when the order of the fractional derivative is $\alpha=1$.

\subsection{Physical Interpretation of the Fractional-Order Continuum Model}
Note that the fractional-order model presented above is based on a continuum mechanics approach. While this is a fundamental approach to mechanics and, in many aspects, analogous to classical and well-established continuum formulations, the model relies on the important hypothesis of a fractional-order kinematics. More specifically, the differential line elements are defined using fractional-order deformation gradients similar to \cite{sumelka2014thermoelasticity,sumelka2016fractional}. This hypothesis results in assuming that the response of a selected point within the solid is affected directly by the response of a collection of points within a characteristic volume, the so-called horizon of nonlocality. Given that the fractional operator is applied directly to the displacement field via the deformation gradient, from a physical standpoint, the formulation accounts for long-range interactions that are proportional to the relative displacement of distant points within the horizon. It follows that, a change in length of an infinitesimal line at the point $\textbf{x}$ between the reference and the current configurations would be affected directly by the response of the points within the nonlocal horizon of $\textbf{x}$. Given the differ-integral nature of space-fractional derivatives, they can be used to capture directly this change in the length of a differential line element at a point $\textbf{x}$ which, in the nonlocal solid, is affected by the response of the points lying in the horizon of $\textbf{x}$. This is indeed a possible formulation of the concept of action-at-a-distance that is often implemented in terms of long-range cohesive forces. In the following, we will show analytically this effect by considering the fractional-order strain.

Using the definition for the Riesz-Caputo fractional derivatives, Eq.~(\ref{infinitesimal_fractional_strain}) can be recast as:
\begin{equation}
\label{eq: comp_step1}
    \tilde{\bm{\epsilon}}(\mathbf{x})={\int_{\textbf{x}-\textbf{l}_A}^{\textbf{x}+\textbf{l}_B}~ \mathcal{A}(\textbf{x},\textbf{s},\textbf{l}_A,\textbf{l}_B,\alpha) \left[\bm{\epsilon}(\textbf{s})\right]~\mathrm{d}\textbf{s}}
\end{equation}
where the kernel $\mathcal{A}(\textbf{x},\textbf{s},\textbf{l}_A,\textbf{l}_B,\alpha)$ is the $\alpha$-order power-law function for the convolution of the classical integer-order strain $\bm{\epsilon}(\textbf{x})$ over the nonlocal horizon $(\textbf{x}-\textbf{l}_A,\textbf{x}+\textbf{l}_B)$. The above equation may be interpreted as if the strain at a point $\textbf{x}$ was given by the weighted-integral of the integer-order strain at generic points ($\bm{\epsilon}(\textbf{s})$) within the domain of influence. Therefore, the nonlocal strain $\tilde{\bm{\epsilon}}(\textbf{x})$ captures the change in the length of a differential line element at the point $\textbf{x}$ in the nonlocal solid, which as mentioned previously, is directly affected by the nonlocal interactions. To show this, the above equation is expressed as:
\begin{equation}
\label{eq: comp_step2}
\begin{split}
    \tilde{\bm{\epsilon}}(\mathbf{x})&=\underbrace{\left({\int_{\mathbf{x}-\mathbf{\delta}}^{\mathbf{x}}~ \mathcal{A}(\mathbf{x},\mathbf{s},\mathbf{l}_A,\mathbf{l}_B,\alpha)\left[\bm{\epsilon}(\mathbf{x})\right]~\mathrm{d}\mathbf{s}}+{\int_{\mathbf{x}}^{\mathbf{x}+\mathbf{\delta}}~ \mathcal{A}(\mathbf{x},\mathbf{s},\mathbf{l}_A,\mathbf{l}_B,\alpha)\left[\bm{\epsilon}(\mathbf{x})\right]~\mathrm{d}\mathbf{s}}\right)}_{{\text{I}}}\\
    &+\underbrace{\left({\int_{\mathbf{x}-\mathbf{l}_A}^{\mathbf{x}-\mathbf{\delta}}~ \mathcal{A}(\mathbf{x},\mathbf{s},\mathbf{l}_A,\mathbf{l}_B,\alpha)\left[\bm{\epsilon}(\mathbf{s})\right]~\mathrm{d}\mathbf{s}}+{\int_{\mathbf{x}+\mathbf{\delta}}^{\mathbf{x}+\mathbf{l}_B}~ \mathcal{A}(\mathbf{x},\mathbf{s},\mathbf{l}_A,\mathbf{l}_B,\alpha)\left[\bm{\epsilon}(\mathbf{s})\right]~\mathrm{d}\mathbf{s}}\right)}_{\text{II: }{\mathcal{\mathbf{R}}(\bm{\epsilon}(\textbf{x}))}}
    \end{split}
\end{equation}
where $\mathbf{\delta} \ll l_A$ and $\mathbf{\delta} \ll l_B$. The above equation can now be expressed as:
\begin{equation}
\label{eq: modified_constitutive}
    \tilde{\bm{\epsilon}}(\mathbf{x})=\underbrace{\gamma~ \bm{\epsilon}(\mathbf{x})}_{\text{I: Local Strain}} + \underbrace{\mathcal{\mathbf{R}}(\bm{\epsilon}(\mathbf{x}))}_{\text{II: Nonlocal Strain}}
\end{equation}
where $\gamma$ is positive constant and the nonlocal contribution of the domains $(\textbf{x}-\textbf{l}_A,\textbf{x}-\mathbf{\delta})$ and $(\textbf{x}+\mathbf{\delta},\textbf{x}+\mathbf{l}_B)$ to the elastic response at point $\textbf{x}$ is accounted in the functional $\mathcal{\mathbf{R}}(\bm{\epsilon}(\textbf{x}))$. For $\mathbf{\delta} \ll l_A$ and $\mathbf{\delta} \ll l_B$, by using the definition of the RC derivative given in Eq.~(\ref{RC_definition}), it can be shown that:
\begin{equation}
    \label{eq: xi_value}
    \gamma=\frac{1}{2}\left[ \left(\frac{\delta}{l_A}\right)^{1-\alpha} + \left(\frac{\delta}{l_B}\right)^{1-\alpha} \right]
\end{equation}
It follows that the constant $\gamma$ is a function of the neighborhood $\mathbf{\delta}$ and length scales $\textbf{l}_A$ and $\textbf{l}_B$, and is strictly $<1$ for $\mathbf{\delta} \ll l_A$ and $\mathbf{\delta} \ll l_B$. It appears from Eq.~(\ref{eq: modified_constitutive}) that the contribution of the nonlocal component $\mathcal{\textbf{R}}(\bm{\epsilon}(\textbf{x}))$ is considered along with the local component $\bm{\epsilon(\textbf{x})}$, as the constant $\gamma<1$. In this sense, we merely note that the fractional-order formulation is analogous to the two-phase nonlocal model developed in \cite{polizzotto2001nonlocal}.
Therefore, from Eq. \eqref{eq: modified_constitutive}, the geometrical definition for the fractional-order strain may be stated as the parameter defined to capture the change of the length of an infinitesimal line in the nonlocal body. This is done by including the effects of the points within the domain of influence through $\mathcal{\textbf{R}}(\bm{\epsilon}(\textbf{x}))$, in addition to the response at the point under study captured by $\bm{\epsilon}(\textbf{x})$. Further, in Appendix 2 we demonstrate the equivalent interpretation of the fractional-order continuum approach and Eringen's model of integral nonlocality.
Additional considerations on the physical meaning of fractional-order models can be found in \cite{di2008long,carpinteri2011fractional,sumelka2016geometrical}.

\subsection{Constitutive Relations of a Nonlocal Euler-Bernoulli Beam}
\label{sec:Constitutive Relations for Nonlocal Beam}
\begin{figure}[ht!]
    \centering
    \includegraphics[width=0.7\textwidth]{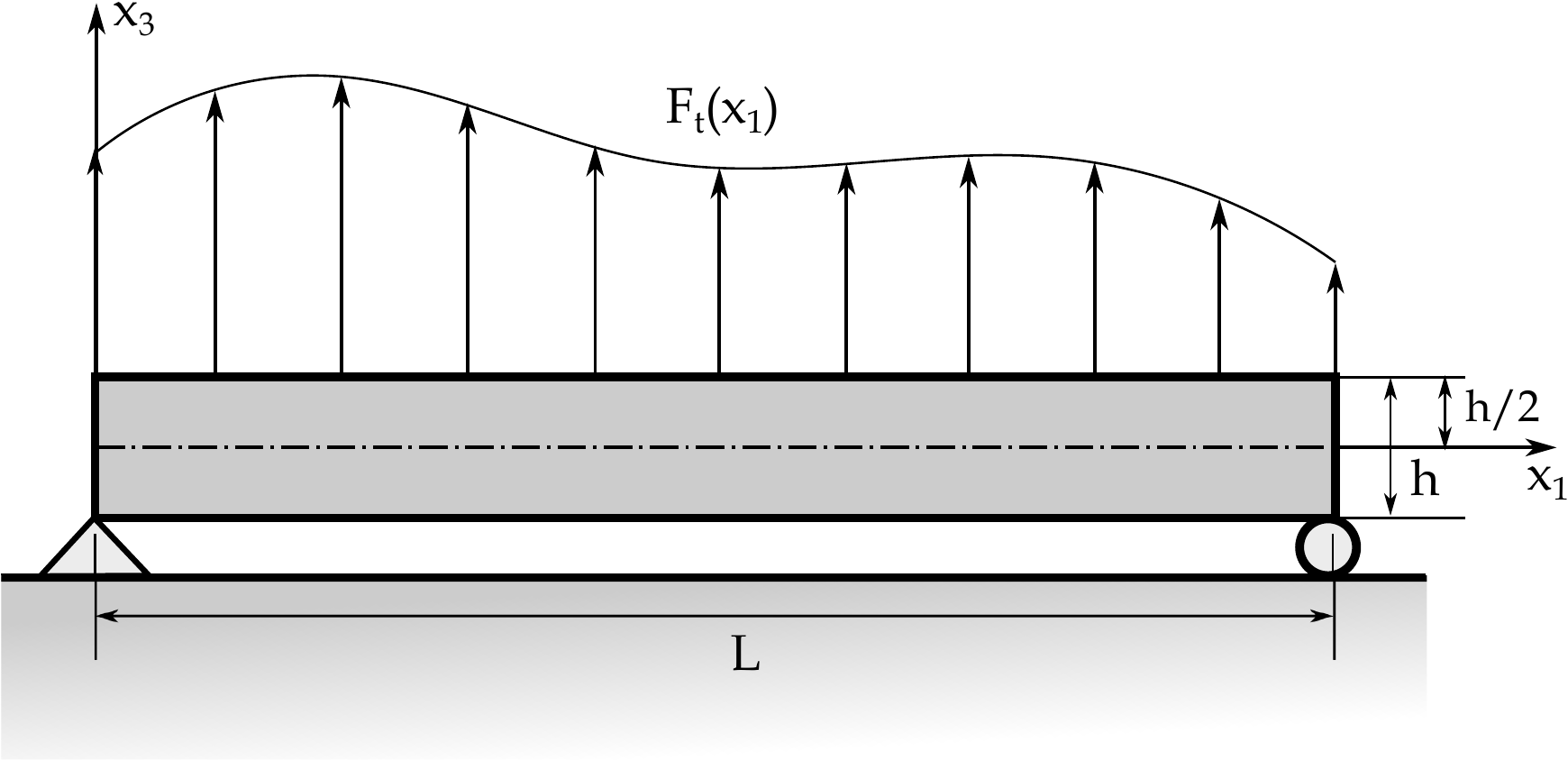}
    \caption{Schematic of a simply-supported linear elastic beam subject to a distributed transverse load $F_t(x_1)$.}
    \label{fig: schematic_beam}
\end{figure}

We use the fractional-order continuum formulation presented above to develop a fractional-order analogue of the constitutive model of a Euler-Bernoulli beam theory for nonlocal slender beams. A schematic of the beam along with the chosen Cartesian reference frame is illustrated in Fig.~(\ref{fig: schematic_beam}). As indicated in the figure, the Cartesian coordinates are chosen such that $x_1=0$ and $x_1=L$ are the longitudinal ends of the beam, and $x_3=\pm h/2$ are the top and bottom surfaces of the beam with $x_3=0$ being the mid-plane of the beam. 

For the chosen coordinate system, the axial and transverse components of the displacement field, denoted by $u_1(x_1,x_3)$ and $u_3(x_1,x_3)$ at any spatial location $\textbf{x}(x_1,x_3)$, are related to the mid-plane displacements according to the Euler-Bernoulli assumptions:
\begin{subequations}
\label{eq: eb_theory}
    \begin{equation}
    u_1(x_1,x_3)=u_0(x_1)-x_3\frac{\mathrm{d}w_0(x_1)}{\mathrm{d}x_1}
    \end{equation}
    \begin{equation}
    u_3(x_1,x_3)=w_0(x_1)
    \end{equation}
\end{subequations}
where $u_0(x_1)$ and $w_0(x_1)$ are the axial and transverse displacements of a point $\textbf{x}(x_1,0)$ on the mid-plane. For the above displacement field, the axial strain $(\tilde{\epsilon}_{11})$ is evaluated using Eq.~(\ref{infinitesimal_fractional_strain}) as:
\begin{equation}
\label{eq: strain_ebt}
    \tilde{\epsilon}_{11}(x_1,x_3)=D_{x_1}^\alpha u_0(x_1)-x_3D_{x_1}^\alpha \frac{\mathrm{d}w_0(x_1)}{\mathrm{d}x_1}
\end{equation}
Note that the above constitutive relations for the beam are similar to \cite{sumelka2015fractional}, although our formulation uses a different definition of the RC derivative, as highlighted in \S\ref{sec:nonlocal_continuum_formulation}. Further, the governing relations in \cite{sumelka2015fractional} have been developed following Newton's approach of force and moment equilibrium rooted in integer-order mathematical models, while, in our study, the governing equations are derived using variational principles. Using the definition of the RC fractional derivative given in Eq.~\eqref{RC_definition}, the fractional derivatives of the axial displacement and the rotation at the mid-plane used in the above Eq.~(\ref{eq: strain_ebt}) are obtained as:
\begin{equation}
\label{eq: rc_disp}
    D_{x_1}^{\alpha} \phi(x_1) =\frac{1}{2}\Gamma(2-\alpha)\left[l_A^{\alpha-1}~\left(~{}^{C}_{x_{A_1}}D^{\alpha}_{x_1} \phi(x_1)\right)-l_B^{\alpha-1}~\left({}^{C}_{x_1}D^{\alpha}_{x_{B_1}} \phi(x_1)\right)\right]
\end{equation}
where $\phi$ is either $u_0$ or $\frac{dw_0}{dx_1}$, and $\textbf{x}_A(x_{A_1},0)$ and $\textbf{x}_B(x_{B_1},0)$ are the terminals of the left- and right-handed Caputo derivatives in Eq.~(\ref{eq: rc_disp}), which coincide with the terminals of the horizon of nonlocality of $\textbf{x}(x_1,0)$ along $\hat{\textbf{x}}_1$. $l_A=x_1-x_{A_1}$ and $l_B=x_{B_1}-x_1$ are the length scales to the left and right of the point $\textbf{x}(x_1,0)$, respectively, along the direction $\hat{\textbf{x}}_1$.

For the Euler-Bernoulli displacement field given in Eq.~\eqref{eq: eb_theory}, using the definition for the nonlocal strain in Eq.~\eqref{infinitesimal_fractional_strain}, a non-zero transverse shear strain would be obtained. However, for the slender beam the rigidity to transverse shear forces is much higher when compared to the bending rigidity. Hence the contribution of the transverse shear deformation towards the deformation energy of the solid can be neglected.
Using the relations between the Lam\'e constants and the elastic modulus of an isotropic solid, the axial stress $\tilde{\sigma}_{11}$ is obtained in terms of the fractional-order axial normal strain using Eq.~(\ref{stress_equation}) as:
\begin{equation}\label{eq: constt_axial}
    \tilde{\sigma}_{11}(x_1,x_3)=E\tilde{\epsilon}_{11}(x_1,x_3)
\end{equation}
where $E$ is the elastic modulus of the solid. Employing the above constitutive relations, the deformation energy $\mathcal{U}$ of the nonlocal elastic beam of volume $\Omega$ is obtained as:
\begin{equation}\label{eq: def_energy}
    \mathcal{U}=\frac{1}{2}\int_{\Omega} \tilde{\sigma}_{11}(x_1,x_3)~\tilde{\epsilon}_{11}(x_1,x_3)~\mathrm{d}V
\end{equation}
Using the above constitutive model for the fractional-order nonlocal elastic beam, we derive the governing differential equations in their strong form and the associated boundary conditions by imposing optimality conditions on the total nonlocal potential energy functional given as:
\begin{equation}
    \label{total_energy_functional}
    \Pi[\textbf{u}(\textbf{x})]=\mathcal{U}-\int_{L}{F_a}(x_1)~u_0(x_1)\mathrm{d}x_1 - \int_{L}{F_t}(x_1)~w_0(x_1)\mathrm{d}x_1
\end{equation}
where the additional integrals in the above expression correspond to the work done by the axial force  $F_a(x_1)$ and the transverse force $F_t(x_1)$, which are applied externally on the plane perpendicular to mid-plane of the beam at $\textbf{x}(x_1,0)$. In the above system, we have also assumed that there are no body forces applied.

\subsection{Governing equations}
\label{sec:Governing equations}
The quasi-static nonlocal elastic response of the beam modeled by the fractional-order continuum model described above is obtained by solving the following system of equations:
\begin{subequations}
\label{all_governing_equations}
\begin{equation}\label{eq: axial_gde}
        {}^{R-RL}_{x_1-l_B}D^{\alpha}_{x_1+l_A}N(x_1)+F_a(x_1)=0~~\forall~~x_1\in(0,L)
\end{equation}
\begin{equation}\label{eq: transverse_gde}
        \frac{\mathrm{d}}{\mathrm{d}x_1}\big[{}^{R-RL}_{x_1-l_B}D^{\alpha}_{x_1+l_A}M(x_1)\big]+F_t(x_1)=0~~\forall~~x_1\in(0,L)
\end{equation}
and subject to the boundary conditions:
\begin{equation}\label{eq: axial_bcs}
    N(x_1)=0~~\text{or}~~\delta u_0(x_1)=0~~\text{at}~~x_1\in\{0,L\}
\end{equation}
\begin{equation}\label{eq: transverse_moment_bcs}
    M(x_1)=0~~\text{or}~~\delta\Bigg[\frac{\mathrm{d}w_0(x_1)}{\mathrm{d}x_1}\Bigg]=0~~\text{at}~~x_1\in\{0,L\}
\end{equation}
\begin{equation}\label{eq: transverse_force_bcs}
    \frac{\mathrm{d}M(x_1)}{\mathrm{d}x_1}=0~~\text{or}~~\delta w_0(x_1)=0~~\text{at}~~x_1\in\{0,L\}
\end{equation}
where $N(x_1)$ and $M(x_1)$ are the axial and bending stress resultants, respectively, given as:
\begin{equation}
\label{eq: stress_resultant_def_1}
    N(x_1)=\int_{-b/2}^{b/2}\int_{-h/2}^{h/2}\tilde{\sigma}_{11}~\mathrm{d}x_3~\mathrm{d}x_2
\end{equation}
\begin{equation}
\label{eq: stress_resultant_def_2}
    M(x_1)=\int_{-b/2}^{b/2}\int_{-h/2}^{h/2}x_3~\tilde{\sigma}_{11}~\mathrm{d}x_3~\mathrm{d}x_2
\end{equation}
\end{subequations}
where $b$ denotes the width of the beam. In the Eqs.~(\ref{eq: axial_gde},\ref{eq: transverse_gde}), ${}^{R-RL}_{x_1-l_B}D^{\alpha}_{x_1+l_A}(\cdot)$ is a Riesz-type Riemann-Liouville (R-RL) derivative of order $\alpha$ defined analogous to the RC derivative in Eq.~\eqref{RC_definition} as:
\begin{equation}\label{eq: riesz_rl_der}
    {}^{R-RL}_{x_1-l_B}D^{\alpha}_{x_1+l_A}(\cdot)=\frac{1}{2}\Gamma(2-\alpha)\left[l_B^{\alpha-1}\left({}^{RL}_{x_1-l_B}D^{\alpha}_{x_1}(\cdot)\right)-l_A^{\alpha-1}\left({}^{RL}_{x_1}D^{\alpha}_{x_1+l_A}(\cdot)\right)\right]
\end{equation}
where ${}^{RL}_{x_1-l_B}D^{\alpha}_{x_1}$ and ${}^{RL}_{x_1}D^{\alpha}_{x_1+l_A}$ are the left- and right-Riemann Liouville fractional derivatives of order $\alpha$. The R-RL derivatives follow from the process of obtaining the governing equations for the beam via variational principles. In this process, it is shown that the adjoint operator for the Riesz-Caputo (RC) fractional derivative, present in the definition of the fractional-order strain, is the R-RL fractional derivative defined in Eq.~(\ref{eq: riesz_rl_der}). The detailed derivations have been provided in Theorem \#3. However, we emphasize that the governing equations contain both the RC and the R-RL fractional derivatives. Note that the governing equations given in Eqs.~(\ref{eq: axial_gde},\ref{eq: transverse_gde}) have been expressed in terms of the stress and moment resultants defined in Eqs.~(\ref{eq: stress_resultant_def_1},\ref{eq: stress_resultant_def_2}). Expressing the governing equations in terms of the displacement field variables, by using the constitutive stress-strain relations of the beam material along with the definitions of the stress and moment resultants, it is immediate that the governing equations also consist of RC derivatives (see Eq.~(\ref{axial_operator_equation})).

Note that the governing equations for the axial and transverse displacements are uncoupled, similar to what is seen in the classical (local) Euler-Bernoulli beam formulation. Further, as expected, the classical Euler-Bernoulli beam governing equations and boundary conditions are recovered for $\alpha=1$.
Note that the solution of the above equations yield mid-plane axial and transverse displacements $u_0(x_1)$ and $w_0(x_1)$. The entire displacement field of the beam can then be obtained using Eq.~(\ref{eq: eb_theory}). We emphasize that since the mapping, $\textbf{x}=\bm{\Psi}(\textbf{X})$, in Eq.~(\ref{motion_description}) is continuous and invertible, the displacement field $\textbf{u}(\textbf{x})$ also belongs to a class $\psi$ of all kinematically admissible displacement fields such that every $\textbf{u}(\textbf{x})\in\psi$ is continuous and satisfies the boundary conditions in the Eqs.~(\ref{eq: axial_bcs},\ref{eq: transverse_moment_bcs}). 
With this condition on the admissible displacement fields we prove the following:\\

\noindent\textbf{Theorem \#1.} \textit{The set of linear operators describing the governing differential equations~(\ref{eq: axial_gde}-\ref{eq: stress_resultant_def_2}) of the beam are self-adjoint.}\\

\noindent\textbf{Proof.} First, we present the proof for the self-adjointness of the operator of the governing equation representing axial motion. Using the expression for the stress resultant in Eq.~(\ref{eq: stress_resultant_def_1}), the governing equation in the axial direction (Eq.~(\ref{eq: axial_gde})) is recast into the following:
\begin{subequations}
\label{axial_operator_equation}
\begin{equation}
    {}^{R-RL}_{x_1-l_B}D^{\alpha}_{x_1+l_A}\left[D_{x_1}^{\alpha}u_0(x_1)\right]+\frac{F_a(x_1)}{Ebh}=0~~\forall~~x_1\in(0,L)
\end{equation}
and subject to the boundary conditions:
\begin{equation}
    \frac{\mathrm{d}u_0(x_1)}{\mathrm{d}x_1}=0~~\text{or}~~\delta u_0(x_1)=0~~\text{at}~~x_1\in\{0,L\}
\end{equation}
\end{subequations}
The operator pertaining to the above fractional-order differential equation is chosen as:
\begin{equation}
\label{operator_definition}
    \tilde{\mathbb{L}}(\cdot)={}^{R-RL}_{x_1-l_B}D^{\alpha}_{x_1+l_A}\left[D_{x_1}^{\alpha}(\cdot)\right]
\end{equation}
Note that the fractional operator $\tilde{\mathbb{L}}(\cdot)$ is linear in nature \cite{podlubny1998fractional}. We consider the inner-product $\langle\tilde{\mathbb{L}}(u_0),v_0\rangle$ such that $u_0$ and $v_0$ satisfy the boundary conditions given in Eq.~(\ref{axial_operator_equation}b):
\begin{equation}
\label{self_adjointness_step_1}
    \langle\tilde{\mathbb{L}}(u_0),v_0\rangle=\int_0^L v_0(x_1)~ {}^{R-RL}_{x_1-l_B}D^{\alpha}_{x_1+l_A}\left[D_{x_1}^{\alpha}u_0(x_1)\right]\mathrm{d}x_1
\end{equation}
Using the definition of the Riesz-type Riemann-Liouville fractional derivative given in Eq.~(\ref{eq: riesz_rl_der}) the above integration is expressed as:
\begin{equation}
\label{self_adjointness_step_2}
    \langle\tilde{\mathbb{L}}(u_0),v_0\rangle=\frac{1}{2}\Gamma(2-\alpha)\int_0^L v_0(x_1)~\left[l_B^{\alpha-1}\frac{\mathrm{d}}{\mathrm{d}x_1}\int_{x_1-l_B}^{x_1}\frac{D_{s_1}^{\alpha}u_0(s_1)}{(x_1-s_1)^{\alpha}}~\mathrm{d}s_1+l_A^{\alpha-1}\frac{\mathrm{d}}{\mathrm{d}x_1}\int_{x_1}^{x_1+l_A}\frac{D_{s_1}^{\alpha}u_0(s_1)}{(s_1-x_1)^{\alpha}}~\mathrm{d}s_1\right]\mathrm{d}x_1
\end{equation}
We further evaluate the above integrals using integration by parts to obtain the following:
\begin{equation}
\label{self_adjointness_step_3}
    \langle\tilde{\mathbb{L}}(u_0),v_0\rangle=\frac{1}{2}\Gamma(2-\alpha)\int_0^L \frac{\mathrm{d}v_0(x_1)}{\mathrm{d}x_1}~\left[l_B^{\alpha-1}\int_{x_1-l_B}^{x_1}\frac{D_{s_1}^{\alpha}u_0(s_1)}{(x_1-s_1)^{\alpha}}\mathrm{d}s_1+l_A^{\alpha-1}\int_{x_1}^{x_1+l_A}\frac{D_{s_1}^{\alpha}u_0(s_1)}{(s_1-x_1)^{\alpha}}\mathrm{d}s_1\right]\mathrm{d}x_1 +\left.v_0\frac{\mathrm{d}u_0}{\mathrm{d}x_1}\right\vert_0^L
\end{equation}
We again highlight that the boundary conditions obtained during the simplification of Eq.~(\ref{self_adjointness_step_2}) to Eq.~(\ref{self_adjointness_step_3}) are integer-order due to the definition of the RC derivative (see Eqs.~(\ref{sq10}-\ref{sq13})). We exchange the order of integration in the above integrals and further, use the boundary conditions in Eq.~(\ref{axial_operator_equation}b) to obtain the following expression:
\begin{equation}
    \langle\tilde{\mathbb{L}}(u_0),v_0\rangle=\frac{1}{2}\Gamma(2-\alpha)\int_0^L D_{s_1}^{\alpha}u_0(s_1)~\left[l_B^{\alpha-1}\int_{s_1}^{s_1+l_B}\frac{D_{x_1}^{1}v_0(x_1)}{(x_1-s_1)^{\alpha}}\mathrm{d}x_1+l_A^{\alpha-1}\int_{s_1-l_A}^{s_1}\frac{D_{x_1}^{1}v_0(x_1)}{(s_1-x_1)^{\alpha}}\mathrm{d}x_1\right]\mathrm{d}s_1
\end{equation}
Using the definition of the RC derivative given in Eq.~(\ref{RC_definition}), the above integral is simplified as:
\begin{equation}
\label{adjoint_step_1}
    \langle\tilde{\mathbb{L}}(u_0),v_0\rangle=\int_0^L D_{x_1}^{\alpha}[u_0(x_1)]~D_{x_1}^{\alpha}[v_0(x_1)]~\mathrm{d}x_1
\end{equation}
By exploiting the symmetry in the above expression, we can write the following:
\begin{equation}
\label{adjoint_step_2}
    \langle u_0,\tilde{\mathbb{L}}(v_0)\rangle=\int_0^{L} D_{x_1}^{\alpha}[u_0(x_1)]~D_{x_1}^{\alpha}[v_0(x_1)]~\mathrm{d}x_1
\end{equation}
By comparing the Eq.~(\ref{adjoint_step_1}) and Eq.~(\ref{adjoint_step_2}), it is clear that the operator $\tilde{\mathbb{L}}(\cdot)$ is self-adjoint. By following the steps through Eqs.~(\ref{axial_operator_equation}-\ref{adjoint_step_2}), it can be similarly shown that the operator describing the transverse governing equation of the beam is also self-adjoint in nature. For the sake of brevity, we skip the proof here. This establishes the claim in Theorem \#1.

It immediately follows from the Theorem \#1 that the system is positive definite. This can be easily verified by considering $\langle\tilde{\mathbb{L}}(u_0),u_0\rangle$ in the Eq.~(\ref{adjoint_step_1}), which results in a quadratic form within the integral. We emphasize that the self-adjointness and positive definiteness of the system hold independently of boundary conditions. This is particularly exciting, because it is established in the literature that it is not possible to define a self-adjoint quadratic potential energy for the classical integral approach to nonlocal elasticity \cite{reddy2007nonlocal,challamel2014nonconservativeness}. Thus, the fractional-order modeling of nonlocality presents us with a way to model nonlocality while ensuring a self-adjoint positive definite system. We will show in \S\ref{sec:numerical_results_and_discussions} that this characteristic further ensures a consistent softening behaviour of the beam with the increasing level of nonlocality, regardless of the boundary conditions. This is unlike the paradoxical predictions of hardening or absence of nonlocal effects, for certain combinations of boundary conditions, in the nonlocal integral theories for beams presented in \cite{challamel2008small,challamel2014nonconservativeness,khodabakhshi2015unified}.\\

\noindent\textbf{Theorem \#2.} \textit{The displacement field $\textbf{u}(\textbf{x})$ which solves the set of equations and boundary conditions in Eqs.~(\ref{eq: axial_gde}-\ref{eq: stress_resultant_def_2}) (if it exists) is unique in the class $\psi$. Further, the strain and stress fields $\tilde{\bm{\epsilon}}(\textbf{x})$ and $\tilde{\bm{\sigma}}(\textbf{x})$ corresponding to the solution $\textbf{u}(\textbf{x})$ are also unique.}\\

\noindent\textbf{Proof.} The proof of the uniqueness of the solution to Eqs.~(\ref{eq: axial_gde}-\ref{eq: stress_resultant_def_2}) is obtained through the method of contradiction, similar to classical continuum analysis. For this, we assume that there exists two different set of solutions to the Eqs.~(\ref{eq: axial_gde}-\ref{eq: stress_resultant_def_2}) given by $\big\{\textbf{u}^{(1)},\tilde{\bm{\epsilon}}^{(1)},\tilde{\bm{\sigma}}^{(1)}\big\}$ and $\big\{\textbf{u}^{(2)},\tilde{\bm{\epsilon}}^{(2)},\tilde{\bm{\sigma}}^{(2)}\big\}$.
We now consider the difference in the displacement and strain fields. The strain field, $\Delta\tilde{\bm{\epsilon}}=\tilde{\bm{\epsilon}}^{(1)}-\tilde{\bm{\epsilon}}^{(2)}$, is also compatible with zero natural boundary conditions. Similarly, the stress field, $\Delta\tilde{\bm{\sigma}}=\tilde{\bm{\sigma}}^{(1)}-\tilde{\bm{\sigma}}^{(2)}$ caused by $\Delta\tilde{\bm{\epsilon}}$ is in equilibrium with zero axial and transverse forces. Clearly, the difference fields $\{\Delta\tilde{\bm{\epsilon}},\Delta\tilde{\bm{\sigma}}\}$ satisfy the homogeneous part of the governing equations and boundary conditions (\ref{eq: axial_gde},\ref{eq: transverse_gde}). Hence, the principle of virtual work dictates that:
\begin{equation}
\label{uniqueness_step_1}
    \int_\Omega \Delta\tilde{\sigma}_{ij}\Delta\tilde{\epsilon}_{ij}=0
\end{equation}
Using the constitutive relation in Eq.~(\ref{stress_equation}), the integral in Eq.~(\ref{uniqueness_step_1}) simplifies to:
\begin{equation}
\label{uniqueness_step_2}
    \int_\Omega (2\mu\Delta\tilde{\epsilon}_{ij}+\lambda\Delta\tilde{\epsilon}_{kk}\delta_{ij})\Delta\tilde{\epsilon}_{ij}=0
\end{equation}
Given the positive definite nature of the above integral, the equality will hold \textit{if} $\Delta\tilde{\bm{\epsilon}}=0~\forall~\textbf{x}\in\Omega$. Hence $\tilde{\bm{\epsilon}}^{(1)}=\tilde{\bm{\epsilon}}^{(2)}~\forall~\textbf{x}\in\Omega$, and consequently from Eq.~(\ref{stress_equation}), $\tilde{\bm{\sigma}}^{(1)}=\tilde{\bm{\sigma}}^{(2)}~\forall~\textbf{x}\in\Omega$.
Further, we highlight here that given the fractional (nonlocal) nature of the strain it might seem that different displacement fields could produce the same strain at a point. However, note that the attenuation function used within the kernel of the fractional derivative is monotonic in nature. Hence, for a given set of fractional parameters $\{\alpha,l_A,l_B\}$, the monotonicity of the kernel dictates that the displacement field is unique. This establishes the uniqueness of the solution set $\{\textbf{u},\tilde{\bm{\epsilon}},\tilde{\bm{\sigma}}\}$. We emphasize that although we consider only the 2D isotropic beam problem, the argument in Theorem \#2 is applicable to any generalized elasticity problem.\\

\noindent\textbf{Theorem \#3.} \textit{The displacement field $\textbf{u}(\textbf{x})\in\psi$ which solves the Eqs.~(\ref{eq: axial_gde}-\ref{eq: stress_resultant_def_2}) minimizes the total potential energy functional given in Eq.~(\ref{total_energy_functional}) in the class $\psi$, and conversely, the displacement field minimizing the functional in Eq.~(\ref{total_energy_functional}) solves the nonlocal beam governing Eqs.~(\ref{eq: axial_gde}-\ref{eq: stress_resultant_def_2}).}\\

\noindent\textbf{Proof of the first claim.} Let $\textbf{u}^\dagger\in\psi$ be the unique solution to the system of equations (\ref{eq: axial_gde}-\ref{eq: stress_resultant_def_2}). Note that we have implicitly assumed that the solution to Eqs.~(\ref{eq: axial_gde}-\ref{eq: stress_resultant_def_2}) exists. Following the standard process of variational calculus we assume $\textbf{u}=\textbf{u}^\dagger+\delta\textbf{u}$ is another kinematically admissible field such that $\delta\textbf{u}\in\psi^\dagger$. The class $\psi^\dagger$ is similar to the class $\psi$ except for the boundary points $\textbf{x}(x_1,x_3)~\forall~x_1\in\{0,L\}$, where the displacement degrees of freedom $\{u_0,w_0,\mathrm{d}w_0 / \mathrm{d}x_1\}=0$ in context of Eqs.~(\ref{eq: axial_gde}-\ref{eq: stress_resultant_def_2}). We highlight here that all quantities $\square^\dagger$ correspond to the displacement field $\textbf{u}^\dagger$. Under the above conditions Eq.~(\ref{total_energy_functional}) yields:
\begin{equation}
    \label{functional_variation}
    \Pi[\textbf{u}]=\Pi[\textbf{u}^\dagger]+\delta\Pi+\frac{1}{2}\delta^2\Pi
\end{equation}
where $\delta\Pi$ and $\delta^2\Pi$ are the first and second variations of $\Pi$ from $\textbf{u}^\dagger$. Using the Eqs.~(\ref{eq: constt_axial},\ref{eq: def_energy}), the first variation $\delta\Pi$ is obtained as:
\begin{equation}
\delta\Pi=\int_\Omega E~\tilde{\epsilon}^\dagger_{11}(x_1,x_3)~\delta[\tilde{\epsilon}_{11}(x_1,x_3)]~\mathrm{
d}V-\int_L F_a(x_1)~\delta u_0~\mathrm{d}x-\int_L F_t(x_1)~\delta w_0~\mathrm{d}x
\end{equation}
Further, using the Eqs.~(\ref{eq: eb_theory}-\ref{eq: rc_disp}) we simplify $\delta\Pi$ to be:
\begin{equation}
\label{minimization_step_1}
\delta\Pi=\int_{0}^{L}N^\dagger(x_1)~D_{x_1}^{\alpha}[\delta u_0(x_1)] \mathrm{d}x_1-\int_{0}^{L}M^\dagger(x_1)~D_{x_1}^{\alpha}\bigg[\frac{\mathrm{d}\delta w_0(x_1)}{\mathrm{d}x_1}\bigg] \mathrm{d}x_1-\int_{0}^{L} F_a(x_1)\delta u_0\mathrm{d}x-\int_{0}^{L} F_t(x_1)\delta w_0\mathrm{d}x
\end{equation}
Before proceeding to further simplify the first variation $\delta\Pi$, we highlight here that the use of variational calculus for the derivation of Euler-Lagrange equations and the transversality conditions has been carried out in \cite{agrawal2002formulation,agrawal2007fractional} for fixed terminals of the fractional derivatives. 
However, in this study, as described in \S\ref{sec:nonlocal_continuum_formulation}, the terminals of the fractional derivatives (equal to the horizon of nonlocality) are position dependent at points close to material boundaries leading to a truncation of the length scales and asymmetry of the horizon of nonlocality. This possible asymmetry in the terminals of the left- and right-handed derivatives in the RC derivative is taken care while simplifying $\delta\Pi$ in Eq.~(\ref{minimization_step_1}).
Using the definition of RC derivative given in Eq.~\eqref{RC_definition} and standard rules of integration, the first and second integrals in the right-hand side of Eq.~(\ref{minimization_step_1}) are simplified to be:
\begin{subequations}
\label{minimization_step_2}
\begin{equation}
    \int_{0}^{L}N^\dagger(x_1)~D_{x_1}^{\alpha}[\delta u_0(x_1)] \mathrm{d}x_1=-\int_{0}^{L}~\left[{}^{R-RL}_{x_1-l_B}D^{\alpha}_{x_1+l_A}N^\dagger(x_1)\right]\delta u_0(x_1)~ \mathrm{d}x_1+\left.\left[N^\dagger(x_1)~\delta u_0(x_1)\right]\right\vert_{0}^{L}
\end{equation}
\begin{equation}
\begin{split}
    \int_{0}^{L}M^\dagger(x_1)~D_{x_1}^{\alpha}\bigg[\frac{\mathrm{d}\delta w_0(x_1)}{\mathrm{d}x_1}\bigg] \mathrm{d}x_1=\int_{0}^{L}&\frac{\mathrm{d}}{\mathrm{d}x_1}\left[{}^{R-RL}_{x_1-l_B}D^{\alpha}_{x_1+l_A}M^\dagger(x_1)\right]\delta w_0(x_1)~\mathrm{d}x_1+\\&
    \left.\left[M^\dagger(x_1)~\delta \left(\frac{\mathrm{d}w_0(x_1)}{\mathrm{d}x_1}\right)\right]\right\vert_{0}^{L}-\left.\left[\frac{\mathrm{d}M^\dagger(x_1)}{\mathrm{d}x_1}~\delta {w_0(x_1)}\right]\right\vert_{0}^{L}
    \end{split}
\end{equation}
\end{subequations}
The detailed steps leading to the simplifications in Eq.~(\ref{minimization_step_2}) are outlined in the Appendix 3. Using Eqs.~(\ref{all_governing_equations},\ref{minimization_step_1},\ref{minimization_step_2}) it can be shown that $\delta\Pi=0$. Additionally, the second variation $\delta^2\Pi$ is given as:
\begin{equation}
    \label{second_variation_positivity}
    \delta^2\Pi=\int_\Omega E~\delta[\tilde{\epsilon}_{11}(x_1,x_3)]~\delta[\tilde{\epsilon}_{11}(x_1,x_3)]\mathrm{d}V
\end{equation}
For any nontrivial $\delta\textbf{u}\in\psi^\dagger$ we have from the above equation $\delta^2\Pi>0$. This leads us to the inequality:
\begin{equation}
    \label{minimization_final}
\Pi[\textbf{u}]=\Pi[\textbf{u}^\dagger]+\frac{1}{2}\delta^2\Pi\geq\Pi[\textbf{u}^\dagger]~~\forall~~\textbf{u}\in\psi
\end{equation}
Note that the equality holds \textit{iff} $\textbf{u}=\textbf{u}^\dagger$ $\forall~\textbf{x}\in\Omega$. Clearly as claimed in \textit{Theorem \#3}, the displacement field $\textbf{u}^\dagger$ which solves the system of equations (\ref{all_governing_equations}) minimizes the functional $\Pi$ in the class $\psi$.\\

\noindent\textbf{Proof of the second claim.} Let $\textbf{u}^\dagger$ be the unique solution to the minimization problem: $\bm{\min}\Pi[\textbf{u}]$ such that $\textbf{u}\in\psi$. The minimization implies that for any variation $\delta\textbf{u}\in\psi^\dagger$, $\delta\Pi$ evaluated at $\textbf{u}^\dagger$ must be identically zero. The $\delta\Pi$ is evaluated through Eqs.~(\ref{minimization_step_1},\ref{minimization_step_2}) where $\textbf{u}^\dagger$, in the context of this proof, minimizes the functional $\Pi$. Clearly, the stress field corresponding to the displacement field $\textbf{u}^\dagger$ uniquely satisfies the equilibrium Eqs.~(\ref{all_governing_equations}), and thus the set $\{\textbf{u}^\dagger,\tilde{\bm{\epsilon}}^\dagger,\tilde{\bm{\sigma}}^\dagger\}$ solves the fractional-order Euler-Bernoulli beam equations.

\section{Fractional Finite Element Method (f-FEM)}
\label{sec:FEM}
The fractional-order nonlocal boundary value problem described by the Eqs.~(\ref{eq: axial_gde}-\ref{eq: transverse_force_bcs}) are numerically solved via a nonlocal finite element method. As discussed in \S \ref{sec:Introduction}, the FE formulation developed for the fractional-order boundary value problem builds upon the FE methods developed in the literature for integral models of nonlocal elasticity \cite{pisano2009nonlocal,norouzzadeh2017finite}. However, several modifications are necessary owing to the choice and the behaviour of the attenuation functions used in the definition of the fractional-order derivatives, as well as the nonlocal continuum model adopted in this study.

\subsection{f-FEM formulation}
The f-FEM is formulated starting from a discretized form of the total potential energy functional $\Pi[\textbf{u}(\textbf{x})]$ given in Eq.~(\ref{total_energy_functional}). Adopting the standard formalism of FEM, the domain $\Omega=[0,L]$ (in case of the slender beam, the length along $x_1$ direction) is divided into $N_e$ finite elements denoted as $\Omega_e$ with $e=\{1,..,N_e\}$ such that $\cup^{N_e}_{e=1}\Omega_e=\Omega$ and $\Omega_j\cap\Omega_k=\emptyset~\forall~j\neq k$.
The unknown variables corresponding to the fractional Euler-Bernoulli theory are given as:
\begin{equation}
\label{eq: disp_vector}
    \{\textbf{u}(\textbf{x})\}^T=[u_0(x_1)~~w_0(x_1)]
\end{equation}
The above vector at any point $x_1$ within $\Omega_e$ is evaluated by an interpolation using the corresponding values at nodes of the element $\Omega_e$. From the definitions of the strains given in Eq. \eqref{eq: strain_ebt}, it may be noted that the interpolation of axial displacement may be carried out using a linear Lagrange approximation, while the transverse displacement field would require the cubic Hermite approximation functions. Therefore, the primary variables at $x_1$ may be expressed as:
\begin{equation}
\label{eq: disp_intpl}
    \{u(x_1)\}=[\hat{N}(x_1)]\{X_g\}
\end{equation}
where $\{X_g\}$ is a vector of the global nodal variables, also referred to as the generalized displacement coordinates, and $[\hat{N}(x_1)]$ is a matrix with appropriate shape functions. We have used the hat symbol on the shape function matrix in order to distinguish it from the axial stress resultant $N(x_1)$. From the definition given in Eq. \eqref{infinitesimal_fractional_strain}, the axial strain in the beam is expressed as:
\begin{equation}\label{eq: strain_att}
    \tilde{\epsilon}_{11}(x_1,x_3)=\frac{1}{2}(1-\alpha){l_A}^{\alpha-1}\left[\int_{{x}_1-l_A}^{{x}_1}\frac{\epsilon_{11}^{I}(s_1,x_3)}{(x_1-s_1)^\alpha}\mathrm{d}s_1\right]+\frac{1}{2}(1-\alpha){l_B}^{\alpha-1}\left[\int_{{x}_1}^{{x}_1+l_B}\frac{\epsilon_{11}^{I}(s_1,x_3)}{(s_1-x_1)^\alpha}\mathrm{d}s_1\right]
\end{equation}
where $\epsilon_{11}^{I}(s_1,x_3)$ is the classical strain (first integer-order derivative of the displacement $u_1$) and $s_1$ is a dummy variable in the direction $\hat{\textbf{x}}_1$ used within the definition of the fractional-order derivative. The above expression is recast into the following:
\begin{subequations}
\label{eq: att_fn_def}
\begin{equation}
    \tilde{\epsilon}_{11}(x_1,x_3)=\int_{{x}_1-l_A}^{{x}_1}A_L(x_1,s_1,l_A,\alpha)\epsilon_{11}^{I}(s_1,x_3)\mathrm{d}s_1+\int_{{x}_1}^{{x}_1+l_B}A_R(x_1,s_1,l_B,\alpha)\epsilon_{11}^{I}(s_1,x_3)\mathrm{d}s_1
\end{equation}
where
\begin{equation}
A_L(x_1,s_1,l_A,\alpha)=\frac{1}{2}(1-\alpha){l_A}^{\alpha-1}\frac{1}{|s_1-x_1|^\alpha}
\end{equation}
\begin{equation}
A_R(x_1,s_1,l_B,\alpha)=\frac{1}{2}(1-\alpha){l_B}^{\alpha-1}\frac{1}{|s_1-x_1|^\alpha}
\end{equation}
\end{subequations}
Note that $A_L(x_1,s_1,l_A,\alpha)$ and $A_R(x_1,s_1,l_B,\alpha)$ are functions of the relative distance between the points $x_1$ and $s_1$, and can be interpreted as similar to the attenuation functions used in integral models of nonlocal elasticity. Clearly, the attenuation decays as a power-law in the distance with an exponent equal to the order $\alpha$ of the fractional derivative. 

Given the nonlocal nature, the expression of the strain in Eq.~(\ref{eq: strain_att}) involves a convolution such that all the points within the horizon of nonlocality $(x_1-l_A,x_1+l_B)$ contribute to the strain at the point $x_1 \in \Omega_e$. While obtaining the FE approximation, the nonlocal contributions from the different finite elements $\tilde{\Omega}_e$ in the horizon $(x_1-l_A,x_1+l_B)$ have to be correctly attributed to the corresponding nodes of those elements $(\tilde{\Omega}_e)$ through appropriate connectivity matrices. In order to correctly attribute these nonlocal contribution, it is essential that we express the nonlocal strain in Eq.~(\ref{eq: strain_att}) as an approximation in terms of the global nodal variables $\{X_g\}$. The process to account for these nonlocal contributions to the elements $\tilde{\Omega}_e$ in the horizon of the element $\Omega_e$ is discussed in detail in the following. 

The expression for the integer-order strain $\epsilon_{11}^{I}(x_1,x_3)$ used in the convolution in Eq.~(\ref{eq: strain_att}) is evaluated for the displacement field given in Eq. \eqref{eq: eb_theory} and expressed in terms of the corresponding nodal variables as:
\begin{equation}
\label{eq: int_strain_intpl}
    \epsilon^I_{11}(s_1,x_3)=\{Z(x_3)\}[B(s_1)]\{X_l(s_1)\}
\end{equation}
where $\{Z(x_3)\}=\{1~~-x_3\}$. $\{X_l(s_1)\}$ is a column vector consisting of the generalized displacements of the nodes of the element enclosing the point $s_1$. The subscript $\square_l$ denotes that $X_l(s_1)$ is a vector containing the local nodal variables of a particular element. We highlight that, in this work, we have used a two-noded element. For the two-noded element $\{X_l(s_1)\}$ is given as:
\begin{equation}
\label{dof_matrix}
    \{X_l(s_1)\}^T=\begin{bmatrix}
    u_0^{(1)} & w_0^{(1)} & \frac{d w_0^{(1)}}{\mathrm{d}s_1} & u_0^{(2)} & w_0^{(2)} & \frac{d w_0^{(2)}}{\mathrm{d}s_1} 
    \end{bmatrix}
\end{equation}
where the superscripts $(\cdot)^{(1)}$ and ${(\cdot)}^{(2)}$ indicate the local node numbers for the element containing the point $s_1$. As evident from Eq.~(\ref{dof_matrix}), there are three degrees of freedom per node. $[B(s_1)]$ is given as:
\begin{equation}
    [B(s_1)]=\begin{bmatrix}
    \frac{\mathrm{d}\mathcal{L}_{1}^e}{\mathrm{d}s_1} & 0 & 0 & \frac{\mathrm{d}\mathcal{L}_{2}^e}{\mathrm{d}s_1} & 0 & 0\\
    0 & \frac{\mathrm{d}^2 \mathcal{H}^e_1}{\mathrm{d}s_1^2} & \frac{\mathrm{d}^2 \mathcal{H}^e_2}{\mathrm{d}s_1^2} & 0 & \frac{\mathrm{d}^2 \mathcal{H}^e_3}{\mathrm{d}s_1^2} & \frac{\mathrm{d}^2 \mathcal{H}^e_4}{\mathrm{d}s_1^2}
    \end{bmatrix}
\end{equation}
where $\mathcal{L}_{1}^e$ and $\mathcal{L}_{2}^e$ are the linear Lagrange shape functions, and $\mathcal{H}_k^e$ with $k\in\{1,2,3,4\}$ are the cubic Hermite shape functions for the two nodes of the element enclosing the point $s_1$.

Using the expression for $\epsilon^I_{11}$ derived in Eq.~(\ref{eq: int_strain_intpl}) the nonlocal strain $\tilde{\epsilon}_{11}$ in Eq.~(\ref{eq: att_fn_def}) is obtained as:
\begin{subequations}
\label{strain_intermediate_step_1}
\begin{equation}
    \tilde{\epsilon}_{11}(x_1,x_3)=\int_{{x}_1-l_A}^{{x}_1+l_B}A(x_1,s_1,l_A,l_B,\alpha)\{Z(x_3)\}[B(s_1)]\{X_l(s_1)\}\mathrm{d}s_1
\end{equation}
where
\begin{equation}
A(x_1,s_1,l_A,l_B,\alpha)=\begin{cases}
A_L(x_1,s_1,l_A,\alpha) & ~~ s_1\in{(x_1-l_A,x_1)}\\
A_R(x_1,s_1,l_B,\alpha) & ~~ s_1\in{(x_1,x_1+l_B)}
\end{cases}
\end{equation}
\end{subequations}
Note that the convolution in Eq.~(\ref{strain_intermediate_step_1}) requires the contribution of the generalized nodal displacements ${X_l(s_1)}$ across the horizon of nonlocality such that $s_1\in(x_1-l_A,x_1+l_B)$. ${X_l(s_1)}$ is further linked to the global nodal variables $\{X_g\}$ through a connectivity matrix $[\tilde{\mathcal{C}}]$ in the following manner:
\begin{equation}
\label{convolution_connectivity_equation}
    \{X_l(s_1)\}=[\tilde{\mathcal{C}}(x_1,s_1)]\{X_g\}
\end{equation}
The connectivity matrix $[\tilde{\mathcal{C}}(x_1,s_1)]$ is crafted such that it activates only the contribution of those nodes that lie within the horizon of nonlocality at the point $x_1$, i.e. for all the nodes enclosing $s_1$ such that $s_1\in(x_1-l_A,x_1+l_B)$. More specifically, for a given point $x_1$, $\{X_l(s_1)\}=\{X_l(s_1)\}~\forall~s_1\in(x_1-l_A,x_1+l_B)$ and $\{X_l(s_1)\}=0$ otherwise. Using the connectivity matrix, Eq.~(\ref{strain_intermediate_step_1}) is simplified to:
\begin{subequations}
\label{nonlocal_strain_in_matrices_final}
\begin{equation}
\tilde{\epsilon}_{11}(x_1,x_3)=[Z(x_3)][\tilde{B}(x_1)]\{X_g\}
\end{equation}
where
\begin{equation}
\label{eq: bt_alpha_def}
    [\tilde{B}(x_1)]=\int_{{x}_1-L_A}^{{x}_1+L_B}A(x_1,s_1,l_A,l_B,\alpha)\{B(s_1)\}[\tilde{\mathcal{C}}(x_1,s_1)]\mathrm{d}s_1
\end{equation}
\end{subequations}
The above expression for the nonlocal strain is used along with the constitutive relations in Eq.~(\ref{eq: constt_axial}) to obtain the total deformation energy defined in Eq. \eqref{eq: def_energy} as:
\begin{equation}
    \mathcal{U}=\frac{1}{2}\{X_g\}^T[\tilde{K}]\{X_g\}
\end{equation}
where the nonlocal stiffness matrix $[\tilde{K}]$ is given as:
\begin{subequations}
\label{eq: stiff_mat_def}
\begin{equation}
    [\tilde{K}]=\int_{0}^{L} [\tilde{B}({x}_1)]^T[D][\tilde{B}({x}_1)]\mathrm{d}x_1
\end{equation}
where
\begin{equation}
    [D]=b\int_{-h/2}^{h/2}E [Z_1(x_3)]^T[Z_1(x_3)]\mathrm{d}x_3
\end{equation}
\end{subequations}
Note that the use of the connectivity matrix in Eq.~(\ref{convolution_connectivity_equation}) results in the fact that the global stiffness matrix $\tilde{K}$ is obtained directly in the global form, hence not requiring a separate assembly process for the element stiffness matrices. As discussed in \cite{pisano2009nonlocal}, owing to the existence of cross-stiffness matrices, the assembly of the element stiffness matrices in a nonlocal FEM requires care and it is not as immediate as the case of a local FEM. Although it might appear that this assembly strategy would require the use of larger (global) matrices, we emphasize that simple principles of connectivity are used to avoid the multiplication of large sparse matrices in Eq.~(\ref{nonlocal_strain_in_matrices_final}), similar to what is done in local FEM. We have provided the specific details of the assembly process in \S\ref{numerical_integration_scheme} where the numerical integration procedure is presented in detail.

The work done by the external axial and transverse forces (see Fig.~(\ref{fig: schematic_beam})) on the beam is expressed in the following manner:
\begin{subequations}
\label{discrete_external_energy}
\begin{equation}
    \mathcal{V}=\{F_{ext}\}\{X_g\}\
\end{equation}
where
\begin{equation}
    \{F_{ext}\}=\int_0^L [F_a(x_1)~~F_t(x_1)]~[\hat{N}(x_1)]~\mathrm{d}x_1
\end{equation}
\end{subequations}
The final algebraic equations describing the FE model of the Euler-Bernoulli beam are now obtained by minimizing the total potential energy of the system $\bm{\Pi}=\mathcal{U}-\mathcal{V}$ as:
\begin{equation}
\label{eq: algebraic_eom}
    [\tilde{K}]\{X_g\}=\{F_{ext}\}
\end{equation}
The solution of the above algebraic equations gives the nodal generalized displacement coordinates.

\subsection{Numerical integration and nonlocal matrices}
\label{numerical_integration_scheme}
In this section, we provide details regarding the numerical scheme employed for the integration of the stiffness matrix and of the force vector. As already discussed the evaluation of the nonlocal stiffness matrix involves a convolution of the classical strain $(\epsilon_{11}^I)$ with the attenuation function $(A(x_1,s_1,l_A,l_B,\alpha))$ due to the fractional-order description of the kinematics (see Eq.~(\ref{Fractional_F_net})). Clearly, when compared to the local FEM, an additional integration procedure has to be carried out to account for the nonlocal behavior. 
The nonlocal interactions across a characteristic length (here, the horizon of nonlocality) have already been accounted for numerically in \cite{polizzotto2001nonlocal,pisano2009nonlocal}. However, differently from these studies, the attenuation function $A(x_1,s_1,l_A,l_B,\alpha)$ in the fractional-order model involves a singularity at the point $x_1$ (more specifically, for $x_1=s_1$) due to the nature of the kernel (see Eq.~(\ref{eq: att_fn_def})).
Although several studies have developed Galerkin FEM for fractional-order BVPs, the governing equations used in these studies do not satisfy the key physical constraints including dimensional consistency and frame-invariance of the nonlocal continuum model. As highlighted earlier, the governing equations used in this study are derived from a frame-invariant continuum model that is better suited to model nonlocal solids with finite and asymmetric horizons.

In order to perform the numerical integration of the nonlocal stiffness matrix, we adopt an isoparametric formulation (see \cite{reddy1989introduction}) and introduce a natural coordinate system $\xi_1$. The Jacobian of the transformation $x_1\rightarrow\xi_1$ is given as $J(\xi_1)$. Using the Gauss-Legendre quadrature rule, the nonlocal stiffness matrix $[\tilde{K}]$ is obtained as:
\begin{equation}
\label{eq: int_stiff_mats}
        [\tilde{K}]=\int_{0}^{L}[\tilde{B}(x_1)]^{T}[D][\tilde{B}(x_1)]\mathrm{d}x_1\approx \sum_{i=1}^{N_e}\sum_{j=1}^{N_{GP}}w_j J^{i}[\tilde{B}(\xi_1^{i,j})]^{T}[D][\tilde{B}(\xi_1^{i,j})]
\end{equation}
where $\xi_1^{i,j}$ is the $j-$th Gauss integration point in the $i-$th  element, $w_j$ is the corresponding weight, and $J^i$ is the Jacobian of the transformation for the $i-$th element. $N_{GP}$ is the number of Gauss points used for the numerical integration over the element, and $N_e$ is the total number of elements into which the beam has been discretized. As previously highlighted, $[\tilde{B}(x_1)]$ (equivalently, $[\tilde{B}(\xi_1^{i,j})]$ in the discretized form) involves an additional integration given the presence of the fractional-order (see Eq.~(\ref{eq: bt_alpha_def})). In the following, we provide the details of this numerical integration.

\begin{figure}
    \centering
    \includegraphics[width=0.8\textwidth]{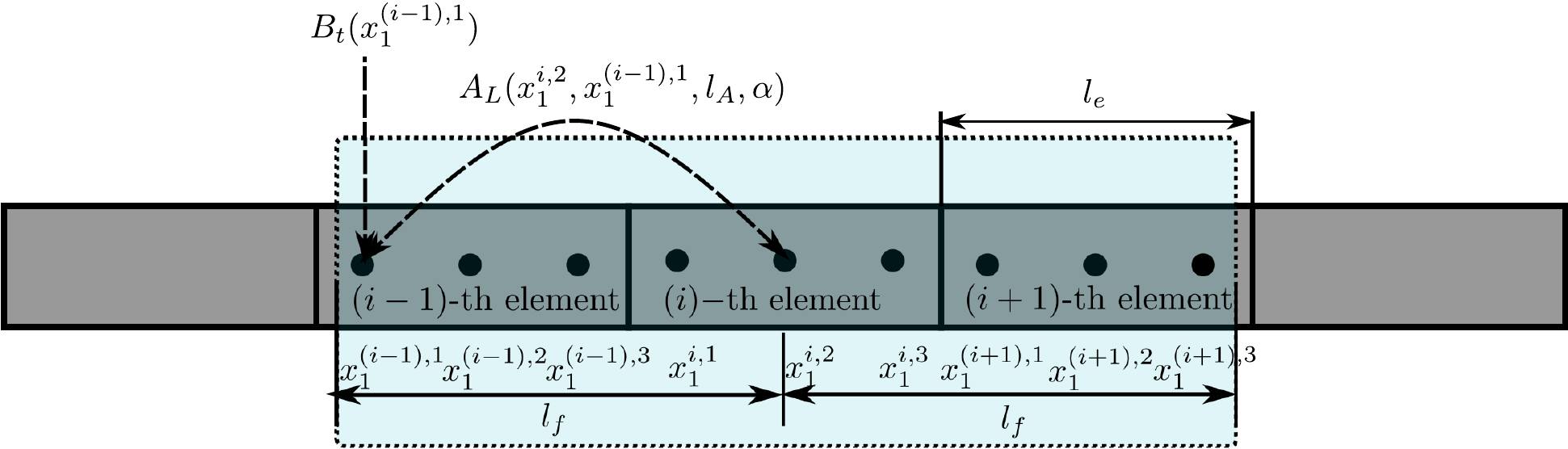}
    \caption{Illustration of the influence zone for the Gauss point $x_1^{i,2}$ and the evaluation of $[\tilde{B}(x_1^{i,2})]$.}
    \label{fig: nlfem_scheme}
\end{figure}

Assume a uniform discretization of the beam such that the length of each element $\Omega^e \in \Omega$ is $l_e$. The integration in $\tilde{B}(x_1)$ has to be performed over the horizon of nonlocality of the point $x_1$ which is given as $(x_1-l_A,x+l_B)$. The number of elements within this horizon to the left $(N_A^{inf})$ and to the right $(N_B^{inf})$ side of $x_1$ are given as $N_A^{inf}=\ceil{l_A/l_e}$ and $N_B^{inf}=\lfloor{l_B/l_e}\rfloor$, respectively. The ceil ($\ceil{\cdot}$) and floor ($\lfloor{\cdot}\rfloor$) functions are used to round the number of elements to the greater integer on the left side and the lower integer on the right side. Using the above formalism, the matrix $[\tilde{B}(x_1)]$ is evaluated in the following manner:
\begin{subequations}
\label{numerical_b_tilde_evaluation}
\begin{equation}
\label{eq: balpha_left_eval}
    \left[\tilde{B}\left(x_1^{(i,j)}\right)\right]=\int_{x_1^{i,j}-l_A}^{x_1^{i,j}}A_L(x^{i,j}_1,s_1,l_A,\alpha)[B(s_1)][\tilde{\mathcal{C}}(x_1^{i,j},s_1)]\mathrm{d}s_1+\int_{x_1^{i,j}}^{x_1^{i,j}+l_B}A_R(x^{i,j}_1,s_1,l_B,\alpha)[B(s_1)][\tilde{\mathcal{C}}(x_1^{i,j},s_1)]\mathrm{d}s_1
\end{equation}
\begin{equation}
\label{eq: balpha_left_eval1}
    \int_{x_1^{i,j}-l_A}^{x_1^{i,j}}A_L(x^{i,j}_1,s_1,l_A,\alpha)[B(s_1)][\tilde{\mathcal{C}}(x_1^{i,j},s_1)]~\mathrm{d}s_1=
    \int_{x_1^{i-N_A^{inf}}}^{x_1^{i-{N_A^{inf}+1}}}\mathcal{I}_L~ \mathrm{d}s_1+...+\int_{x_1^{i-1}}^{x_1^{i}}\mathcal{I}_L~\mathrm{d}s_1+\int^{x_1^{i,j}}_{x_1^{i}}\mathcal{I}_L~\mathrm{d}s_1
\end{equation}
\begin{equation}
\label{eq: balpha_right_eval}
    \int_{x_1^{i,j}}^{x_1^{i,j}+l_B}A_R(x^{i,j}_1,s_1,l_B,\alpha)[B(s_1)][\tilde{\mathcal{C}}(x_1^{i,j},s_1)]~\mathrm{d}s_1=\int_{x_1^{i,j}}^{x_1^{i+1}}\mathcal{I}_R~\mathrm{d}s_1+\int_{x_1^{i+1}}^{x_1^{i+2}}\mathcal{I}_R~\mathrm{d}s_1+...\int_{x_1^{i+N_B^{inf}-1}}^{x_1^{i+{N_B^{inf}}}}\mathcal{I}_R~\mathrm{d}s_1
\end{equation}
\end{subequations}
where $x_1^{i,j}$ is the global coordinate of the Gauss-point $\xi_1^{i,j}$ within the $i-$th element, the ends of which are given as $[x_1^i,x_1^{i+1}]$. For the sake of brevity, we have denoted the integrands in the above equation in the following manner: $ \mathcal{I}_L=A_L(x^{i,j}_1,s_1,l_A,\alpha)[B(s_1)][\tilde{\mathcal{C}}(x_1^{i,j},s_1)]$, and $\mathcal{I}_R=A_R(x^{i,j}_1,s_1,l_B,\alpha)[B(s_1)][\tilde{\mathcal{C}}(x_1^{i,j},s_1)]$. For the elements and the associated Gauss-points that are close to the boundaries ($x_1=0$ or $x_1=L$), the number of elements in the horizon $N_\square^{inf}~(\square\in\{A,B\})$ are truncated appropriately, in order to account for the asymmetric length scales (see \S\ref{sec:NBM}).
Each integration in Eq.~(\ref{numerical_b_tilde_evaluation}) across the elements $\tilde{\Omega}^e$ in the influence zone is performed using the Gauss-Legendre quadrature rule. We emphasize that the strain-displacement matrix for the RC fractional strain-displacement relations given in Eq.~\eqref{eq: strain_ebt} is expressed as a convolution of the integer-order derivatives using $[B(s_1)]$ as shown in Eq.~\eqref{nonlocal_strain_in_matrices_final}. Therefore, the computation of $[B(s_1)]$ is straightforward and follows from Eq.~(\ref{eq: int_strain_intpl}). The integration in Eq.~(\ref{numerical_b_tilde_evaluation}) is now obtained as:
\begin{subequations}
\label{final_b_tilde_evaluation}
\begin{equation}
\label{eq: balpha_left_sing}
 \int_{x_1^{i,j}-l_A}^{x_1^{i,j}}\mathcal{I}_L~\mathrm{d}s_1=\underbrace{\int_{x_1^{i,j}-l_A}^{x_1^{i-N_A^{inf}}}\mathcal{I}_L~\mathrm{d}s_1+...\int_{x_1^{i-1}}^{x_1^{i}}\mathcal{I}_L~\mathrm{d}s_1}_{\text{Gauss-Legendre Quadrature}}+\underbrace{\int_{x_1^{i+1}}^{x_1^{i,j}}\mathcal{I}_L~\mathrm{d}s_1}_{\text{Singularity at } x_1^{i,j}}
 \end{equation}
\begin{equation}
\label{eq: balpha_right_sing}
   \int_{x_1^{i,j}}^{x_1^{i,j}+l_B}\mathcal{I}_R~\mathrm{d}s_1=\underbrace{\int_{x_1^{i,j}}^{x_1^{i+1}}\mathcal{I}_R~\mathrm{d}s_1}_{\text{Singularity at } x_1^{i,j}}+\underbrace{\int_{x_1^{i+1}}^{x_1^{i+2}}\mathcal{I}_R~\mathrm{d}s_1+...\int_{x_1^{i+N_B^{inf}-1}}^{x_1^{i-{N_B^{inf}}}} \mathcal{I}_R~\mathrm{d}s_1}_{\text{Gauss-Legendre Quadrature}}
\end{equation}
\end{subequations}
As discussed earlier, a singularity occurs in the interval containing the Gauss-point $x_1^{i,j}$ for both the left and right integrals due to the nature of the fractional-order derivative kernel. This is unlike the attenuation functions employed for existing models of nonlocal integral elasticity, where no such singularity exists across the domain \cite{polizzotto2001nonlocal,pisano2009nonlocal}. Several numerical strategies have been formulated in the literature to evaluate integrals with end-point singularities, like the composite quadrature, singularity-removing transformations \cite{takahasi1973quadrature}, graded meshes \cite{schwab1994variable}, and adaptive methods \cite{cools2003algorithm}. However, in this work, the end-point singularity in the integral is circumvented by a simple analytical evaluation of that particular integration. 
The Gauss-Legendre quadrature method used to numerically evaluate the other integrals in Eq.~(\ref{final_b_tilde_evaluation}) without singularities is schematically illustrated in Fig.~(\ref{fig: nlfem_scheme}). This integration can also be carried out using other numerical integration procedures including the fast Gauss transform presented in \cite{benvenuti2006fast}. However, the quadrature method allows a simpler implementation of the nonlocal FE code as discussed in \cite{pisano2009nonlocal}. Here below, we provide the expression of this quadrature based integration for the interaction between the $p-$th element lying in the horizon of nonlocality of the Gauss point $x_1^{i,j}$ and $x_1^{i,j}$ itself: 
\begin{equation}
\label{nonlocal_integral_example}
    \int_{x_1^{p}}^{x_1^{p+1}}A_R(x^{i,j}_1,s_1,l_B,\alpha)[B(s_1)][\tilde{\mathcal{C}}(x_1^{i,j},s_1)]~\mathrm{d}s_1=
    \sum_{k=1}^{N_{GP}}w_k J^p A_R(x^{i,j}_1,{x}_1^{p,k},l_B,\alpha)[B({s}_1^{p,k})][\tilde{\mathcal{C}}(x_1^{i,j},{x}_1^{p,k})]
\end{equation}
where ${x}_1^{p,k}$ is the Cartesian coordinate of the $k-$th Gauss point in the $p-$th element which lies in the horizon of nonlocality of the $i-$th element, $w_k$ is the corresponding weight and $J^p$ is the associated Jacobian for this transformation for the $p-$th element. Note that given the nonlocal nature, the attenuation function has to be evaluated using the absolute Cartesian coordinates. In Eq.~(\ref{nonlocal_integral_example}), we have implicitly assumed that the $p-$th element is on the right-hand horizon of $x_1^{i,j}$. However, the same analogy extends directly to evaluate the nonlocal contribution of any element that lies on the left-hand horizon of the point $x_1^{i,j}$. The above algorithm is used for the computation of $[\tilde{B}(x_1^{i,j})]$ at the Gauss-points necessary for the integration in Eq.~(\ref{eq: int_stiff_mats}) to obtain the stiffness matrix. The resulting stiffness matrix is then used in the algebraic equations of motion given in Eq.~(\ref{eq: algebraic_eom}) to determine the generalized displacement coordinates of the fractional-order nonlocal beam.

Before going forward and using the above f-FEM to present numerical simulations, we discuss the assembly procedure of the nonlocal contributions from the elements $\tilde{\Omega}^e$ in the nonlocal horizon of a given point. 
As emphasized earlier, the matrix $\tilde{B}(x_1)$ is obtained directly in the global form. It is clear from the convolution in Eq.~(\ref{nonlocal_integral_example}) that the connectivity matrix $[\tilde{\mathcal{C}}(x_1^{i,j},{x}_1^{p,k})]$ would ensure the proper assembly of the $\tilde{B}(x_1)$ matrix. However, instead of directly multiplying the large connectivity matrix $\tilde{\mathcal{C}}$, its sparse nature is utilized to simplify the computation in the following manner:
\begin{equation}
\label{B_assembly}
    \tilde{B}(:,3n_p-2:3(n_p+1))=\tilde{B}_l^{(p)}
\end{equation}
where $n_p$ and $n_p+1$ are the global node numbers of the $p-$th element in the horizon of the point $x_1^{i,j}$. The colon symbol $(:)$ indicates a sweep operation. More specifically, in the above operation, both the rows of $\tilde{B}_l^{(p)}$ are assigned to both the rows of $\tilde{B}$, and the columns of $\tilde{B}_l^{(p)}$ are added to the column numbers ranging from $3n_p-2$ to $3(n_p+1)$ of the $\tilde{B}$ matrix. $\tilde{B}_l^{(p)}$ is a local matrix corresponding to the $p-$th element in the horizon of the point $x_1^{i,j}$ and is given as:
\begin{equation}
\label{B_local}
    \tilde{B}_l^{(p)}=\int_{x_1^{p}}^{x_1^{p+1}}A_R(x^{i,j}_1,s_1,l_B,\alpha)[B(s_1)]~\mathrm{d}s_1
\end{equation}
Note that the connectivity matrix $[\tilde{\mathcal{C}}(x_1^{i,j},s_1)]$ in Eqs.~(\ref{nonlocal_integral_example}) is absent in the Eq.~(\ref{B_local}) rendering $\tilde{B}_l^{(p)}$ local.
Another remark regarding the assembly procedure pertains to the fact that the Eq.~(\ref{B_assembly}) is valid under the assumptions made in this study of two nodes per element and three degrees of freedom per node. For a general case where the number of nodes in a element is $n_k$ and the number of degrees of freedom per node is $n_d$ the assembly procedure in Eq.~(\ref{B_assembly}) can be modified as:
\begin{equation}
\label{general_B_assembly}
    \tilde{B}(:,n_d n_p-n_d+1:n_d(n_p+n_k-1))=\tilde{B}_l^{(p)}
\end{equation}
where $n_p$ and $n_p+n_k-1$ are the global node numbers of the terminals of the $p-$th element.
The above strategy ensures that the nonlocal contributions from the elements in the horizon of nonlocality of a given Gauss point $x_1^{i,j}$ are correctly attributed to their nodes in the global matrix $\tilde{B}$. Thus, the stiffness matrix $\tilde{K}$ is obtained from Eq.~(\ref{eq: int_stiff_mats}) in the global form requiring no further assembly of the local and nonlocal stiffness contribution matrices.


\section{Numerical Results and Discussion}
\label{sec:numerical_results_and_discussions}
The f-FEM developed above is used to analyze the fractional-order nonlocal Euler-Bernoulli beam under various loading conditions. Further, the effect of the fractional order $(\alpha)$, and of the length scales $(l_A,l_B)$ on the static response of the beam is analyzed here. The aspect ratio of the beam used in this study is fixed at $L/h=100$, in order to satisfy the slender assumption for the Euler-Bernoulli beam displacement theory. The width of the beam is always maintained equal to $h$. The elastic modulus for the isotropic solid is $E=30$ GPa. Note that although the above presented f-FEM is capable of simulating the static response of the fractional-order nonlocal beam under both axial and transverse loads, here we only present the results for the transverse response. Further, in the following, we have assumed a symmetric horizon of nonlocality for points sufficiently inside the domain of the beam, i.e., $l_A=l_B=l_f$. The horizon lengths on the appropriate side of a point are truncated as the point approaches the external boundary, as discussed in \S\ref{sec:NBM}. Before presenting the results of the static response of the fractional-order nonlocal beam for the various loading conditions, we present the results of the studies carried out for validating the f-FEM and establishing its convergence. 

\subsection{Validation}
We validate the f-FEM by comparing the results of the f-FEM against: (a) Validation \#1: exact solution of a fractional-order beam clamped at both its ends and subject to a spatially varying transverse force; (b) Validation \#2: exact solution of a fractional-order beam simply-supported at both its ends and subject to a spatially varying transverse force; and (c) Validation \#3: changing the kernel of the fractional-order derivative to the exponential kernel used in nonlocal integral elasticity \cite{fernandez2016bending,norouzzadeh2017finite} and comparing the obtained nonlocal response with the results available in literature. The specific details of three-fold validation strategy are provided in the following.\\

\begin{figure*}[ht!]
    \centering
    \begin{subfigure}[t]{0.5\textwidth}
        \centering
    \includegraphics[width=\textwidth]{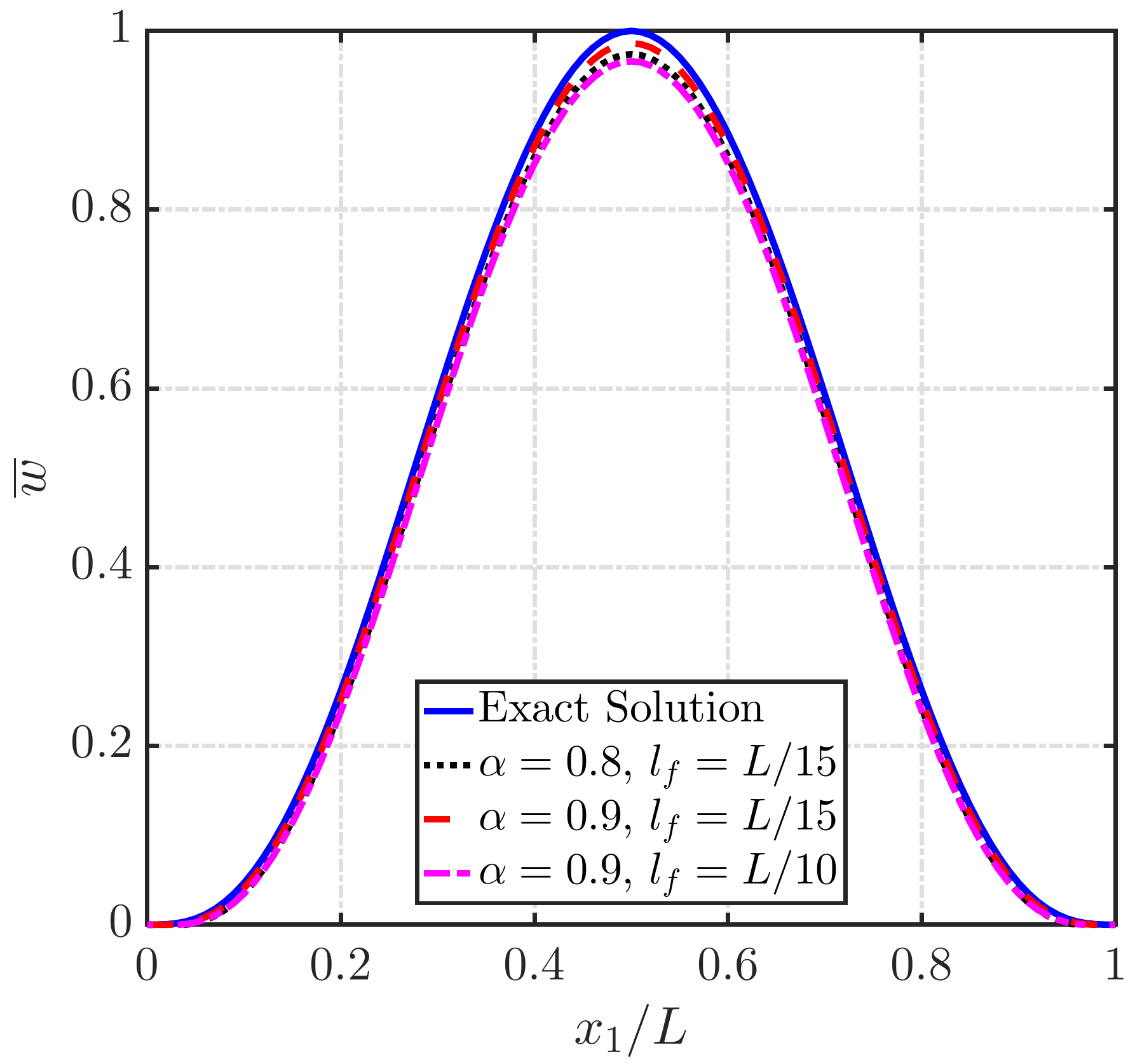}
    \caption{$\overline{w}$.}
    \label{fig: validation_exact_cc_w}
    \end{subfigure}%
    ~ 
    \begin{subfigure}[t]{0.5\textwidth}
        \centering
        \includegraphics[width=\textwidth]{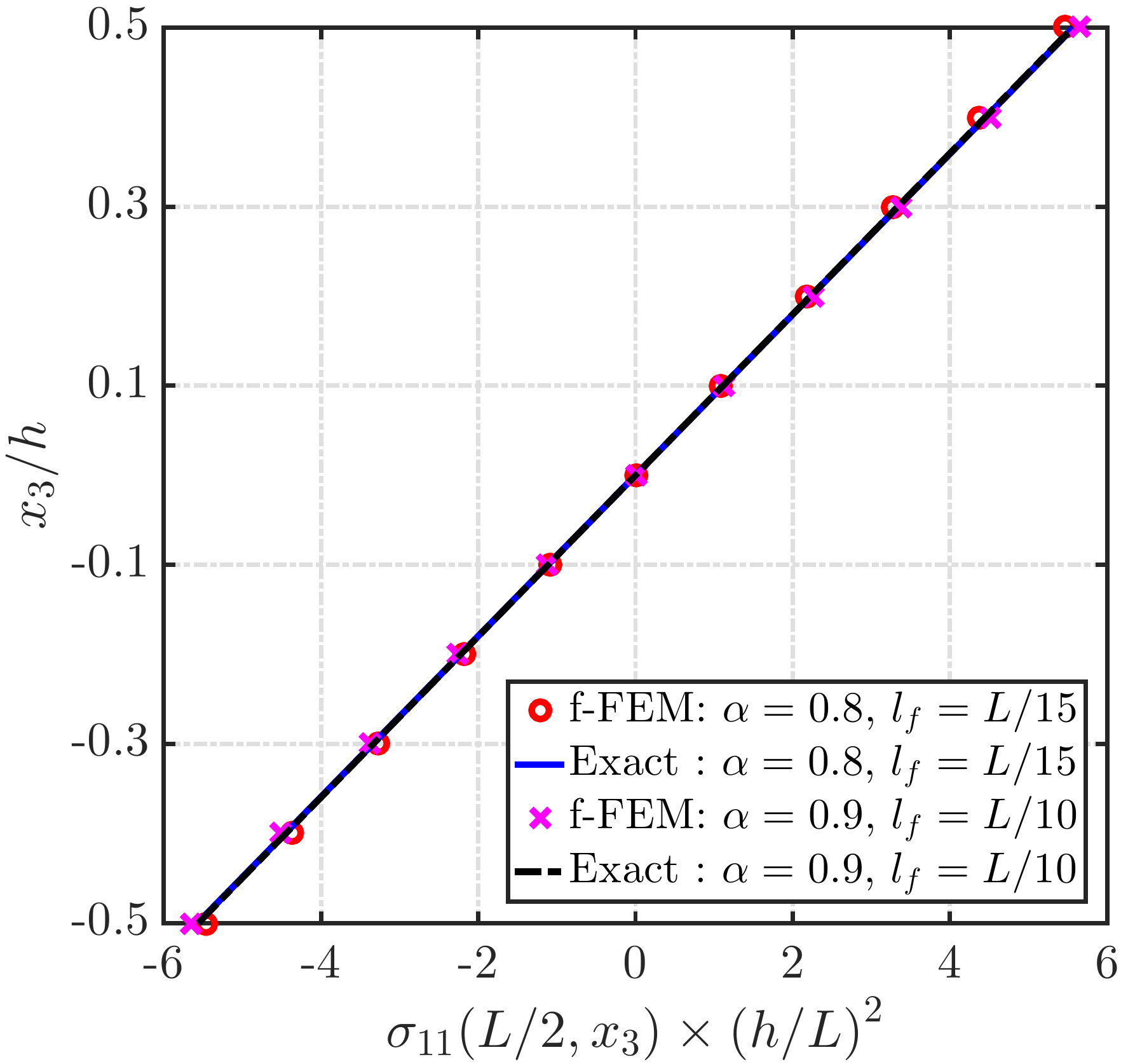}
        \caption{$\overline{\sigma}_{11}$.}
        \label{fig: validation_exact_cc_sig11}
    \end{subfigure}
    \caption{Numerical f-FEM results (a) transverse displacement, and (b) normal axial stress, corresponding to the assumed displacement field in Eq.~\eqref{validation_1_displacement} compared with the exact solutions for the clamped-clamped boundary condition. Note that in validation \#1, we have assumed a transverse displacement field that is independent of $\alpha$ and $l_f$.}
    \label{fig: validation_exact_cc}
\end{figure*}

\noindent \textbf{Validation \#1}: the following transverse displacement of the mid-plane of the beam is assumed:
\begin{equation}
    \label{validation_1_displacement}
    w_0(x_1)=L\left(\frac{x_1}{L}\right)^3\left(1-\frac{x_1}{L}\right)^3
\end{equation}
Note that the above transverse displacement satisfies the boundary conditions for a beam clamped at both ends. By using the strong form of the governing differential equations in Eq.~(\ref{eq: transverse_gde}) the transverse load required for the above displacement response is obtained as:
\begin{equation}
    \label{validation_1_force}
     F_t(x_1)=-\frac{6Eh^3}{L^3}\left[1-5 \left(\frac{x_1}{L}\right)+5 \left(\frac{x_1}{L}\right)^2+10 \left(\frac{l_f}{L}\right)^2\left( \frac{1-\alpha}{3-\alpha}\right)\right]
\end{equation}
We highlight here that the assumed transverse displacement field is independent of $\alpha$ and $l_f$, hence resulting in a forcing function which is dependent on the fractional parameters.
The transverse force distribution is then used within the f-FEM and the numerical approximation for the transverse displacement obtained from the f-FEM is compared against the exact solution given in Eq.~(\ref{validation_1_displacement}) for different combinations of the order $(\alpha)$ and the length scale $(l_f)$. The number of elements used in the discretization of the beam for generating the numerical results is maintained at ten times the ratio of the total length $(L)$ and the length scale $(l_f)$, i.e., $N_e=10(L/l_f)$. We have established in \S\ref{ssec:convergence} that the f-FEM converges for this assumption of $N_e$. The numerical results in terms of the transverse displacement as well as the axial stress are presented in Fig.~(\ref{fig: validation_exact_cc}). The transverse displacement presented in Fig.~(\ref{fig: validation_exact_cc}) is normalized in the following manner:
\begin{equation}
    \overline{w}(x_1)=\frac{64}{L}w_0(x_1)
\end{equation}
As evident from Fig.~(\ref{fig: validation_exact_cc}), the match between the transverse displacement $\overline{w}$ and the axial normal stress $\sigma_{11}$ obtained numerically from the f-FEM and the exact solution is excellent. The error between the numerically obtained f-FEM and the exact solutions is less than $5\%$ in all the cases.\\

\begin{figure*}[ht!]
    \centering
    \begin{subfigure}[t]{0.5\textwidth}
        \centering
    \includegraphics[width=\textwidth]{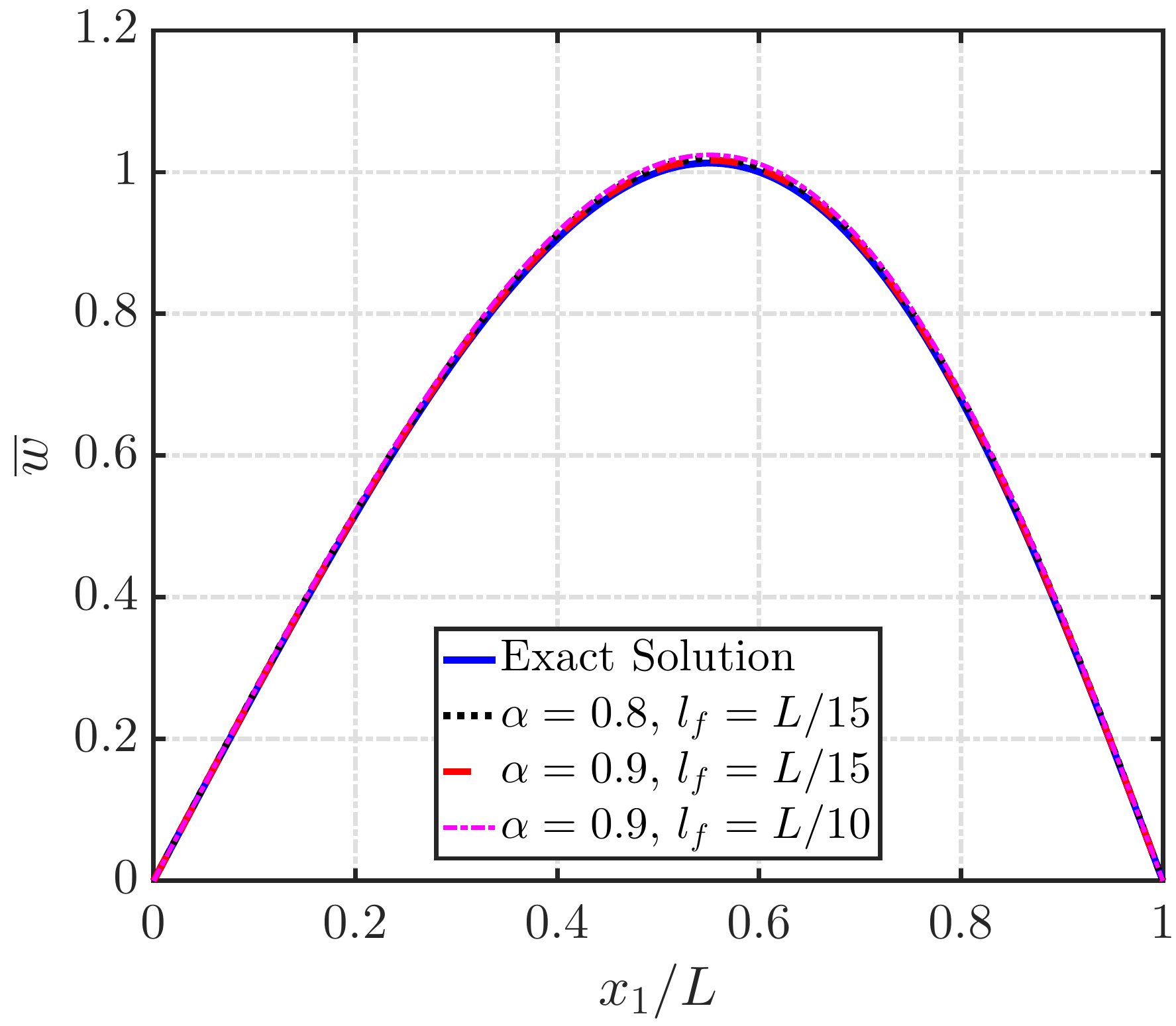}
    \caption{$\overline{w}$.}
    \label{fig: validation_exact_ss_w}
    \end{subfigure}%
    ~ 
    \begin{subfigure}[t]{0.5\textwidth}
        \centering
        \includegraphics[width=\textwidth]{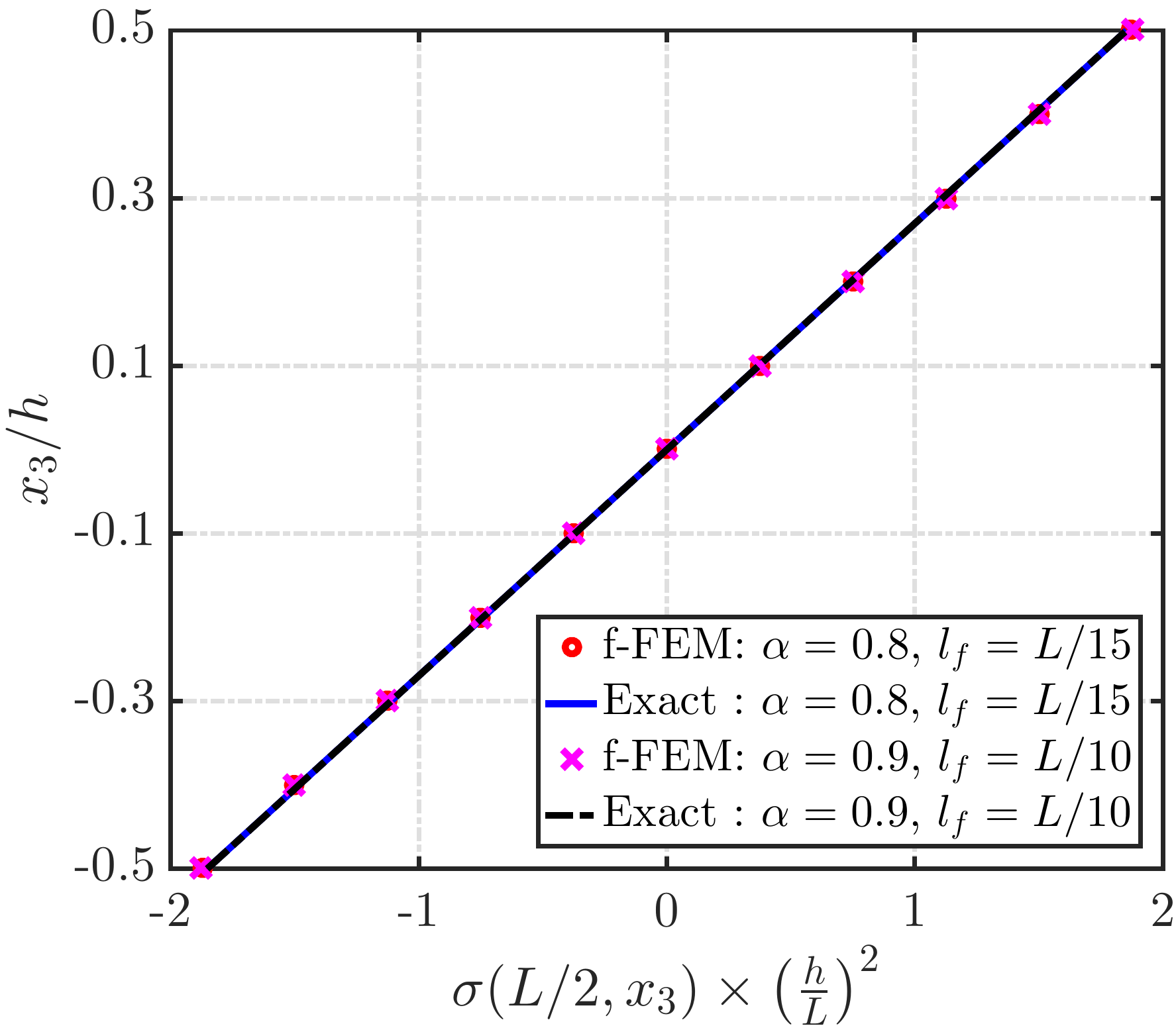}
        \caption{$\overline{\sigma}_{11}$.}
        \label{fig: validation_exact_ss_sig11}
    \end{subfigure}
    \caption{Numerical f-FEM results (a) transverse displacement, and (b) normal axial stress, corresponding to the assumed displacement field in Eq.~\eqref{validation_2_displacement} compared with the exact solutions for the simply-supported boundary condition. Note that in validation \#2, we have assumed a transverse displacement which is independent of $\alpha$ and $l_f$.}
    \label{fig: validation_exact_ss}
\end{figure*}

\noindent \textbf{Validation \#2}: repeating the strategy outlined in validation \#1, the following expression is assumed to be the mid-plane transverse displacement of a simply-supported beam:
\begin{equation}
    \label{validation_2_displacement}
    w_0(x_1)=\frac{L}{100}\left[\left(\frac{x_1}{L}\right)^6+\frac{3}{2}\left(\frac{x_1}{L}\right)^5-4\left(\frac{x_1}{L}\right)^4-2\left(\frac{x_1}{L}\right)^3+\frac{7}{2}\left(\frac{x_1}{L}\right)\right]
\end{equation}
The above transverse displacement satisfies the essential and natural boundary conditions for a beam simply-supported at both ends. Again, by using the strong form of the governing differential equations given in Eq.~(\ref{eq: transverse_gde}) the required transverse load is obtained as:
\begin{equation}
    \label{validation_2_force}
     F_t(x_1)=-\frac{Eh^3}{1200L^3}\left[360\left(\frac{x_1}{L}\right)^2+180\left(\frac{x_1}{L}\right)+720\left(\frac{l_f}{L}\right)^2\bigg(\frac{1-\alpha}{3-\alpha}\bigg)-96\right]
\end{equation}
The above transverse force distribution is used within the f-FEM and the numerical approximation of the transverse displacement is obtained for different combinations of $\alpha$ and $l_f$. The numerical results obtained are presented in Fig.~(\ref{fig: validation_exact_ss}) and compared against the exact result. The number of elements used in the numerical simulation is also maintained at $N_e=10(L/l_f)$. Both the transverse displacement and the axial stress of the mid-plane are compared in Fig.~(\ref{fig: validation_exact_ss}). The transverse displacement is normalized as:
\begin{equation}
    \overline{w}(x_1)=\frac{16}{21L}w_0(x_1)
\end{equation}
As evident from Fig.~(\ref{fig: validation_exact_ss}) the match between the numerically obtained f-FEM results and the exact solution is excellent and the error is less than $3\%$ in all the cases.\\

\noindent \textbf{Validation \#3}: We finally validate the presented f-FEM by simulating a nonlocal beam with an exponential attenuating function commonly used in nonlocal integral elasticity methods.
The constitutive relations for the nonlocal integral elasticity are given as \cite{eringen1972linear}:
\begin{equation}\label{eq: eringen_nonlocal}
    \sigma_{ij}(\mathbf{x})=\int_{\Omega}A(\mathbf{x},\textbf{s},l)~C_{ijkl}~\epsilon_{kl}(\textbf{s})~\mathrm{d}{\mathbf{s}}
\end{equation}
where $A(\mathbf{x},\textbf{s},l)$ is the attenuation function typical of nonlocal integral elasticity models, $C_{ijkl}$ is the constitutive matrix of the solid, and $\epsilon_{kl}$ is the integer-order classical strain. For this model, the total deformation energy necessary to develop the algebraic equations of equilibrium for the weak form is given by:
\begin{equation}
    \Pi^{\text{Eringen}}=\int_{\Omega} \sigma_{ij}(\mathbf{x})~\epsilon_{ij}(\mathbf{x})~\mathrm{d}V
\end{equation}
Note that the strain field in the nonlocal integral model of elasticity is still local in nature, while in the fractional-order formulation, the strain field is nonlocal in nature due to the fractional-order derivatives (see Eq.~(\ref{infinitesimal_fractional_strain})). Therefore, the nonlocal stiffness matrix in Eq.~(\ref{eq: int_stiff_mats}) is modified as:
\label{validation_stiff_mat_def}
\begin{equation}
    [\tilde{K}_v]=\int_{0}^{L} [\tilde{B}({x}_1)]^T[D][{B}({x}_1)]\mathrm{d}x_1
\end{equation}
We highlight that the above stiffness matrix $[\tilde{K}_v]$ involves the multiplication of a nonlocal $[\tilde{B}({x}_1)]$ and a local $[{B}({x}_1)]$ strain-displacement matrix, while the stiffness matrix in the fractional-order nonlocal modeling involves only the nonlocal $[\tilde{B}({x}_1)]$ matrices (see Eq.~(\ref{eq: stiff_mat_def})). Further the attenuation function $A(x_1,s_1,l_f)$ given in Eq. \eqref{eq: bt_alpha_def} for the fractional-order model of nonlocal elasticity is modified to be the exponential order attenuation function given below \cite{fernandez2016bending,norouzzadeh2017finite}:
\begin{equation}\label{eq: eringen_att}
    A(x_1,s_1,l_f)=\frac{1}{2l_f}\exp{-\frac{|x_1-s_1|}{l_f}}
\end{equation}
where the variables retain their definitions. The horizon of nonlocality in the nonlocal elasticity problem modeled in \cite{fernandez2016bending,norouzzadeh2017finite} is assumed to be the entire length of the beam. Thus, another modification involves expanding the horizon of nonlocality beyond the interval $(x_1-l_A,x+l_b)$ to the entire length of the beam while computing $[\tilde{B}(x_1)]$ in Eq.~(\ref{eq: eringen_att}). With these modifications, the numerical results from the f-FEM are compared with the results provided in \cite{fernandez2016bending,norouzzadeh2017finite} for a uniformly distributed transverse load (UDL) on a cantilever beam. For brevity, we only report the comparisons between the maximum displacements which are obtained at the tip of the cantilever beam. The results are summarized in Table \ref{tab: eringen_validation}. For this study, we have normalized the transverse displacement similar to \cite{norouzzadeh2017finite}. Note that the results presented in \cite{norouzzadeh2017finite,fernandez2016bending} have been obtained from the nonlocal analogue of the Timoshenko beam theory which is a shear-deformable theory. In the present work, we have adopted a Euler-Bernoulli theory. As already established in classical theories, shear-deformable beam models exhibit lower stiffness when compared to the Euler-Bernoulli ones. This observation is consistent with the results in our study that show lower transverse deformation than that of \cite{norouzzadeh2017finite}. Nevertheless, given that these results are obtained for slender beams, the match is excellent and an error $\leq 1\%$ is obtained in all the cases (see Table \ref{tab: eringen_validation}).

\begin{table}[hb!]
    \centering
        \caption{Nonlocal elastic response for a cantilever beam having an aspect ratio $L/h=25$ subjected to UDL from the f-FEM validated with nonlocal integral elasticity results from \cite{fernandez2016bending,norouzzadeh2017finite} ($\overline{w}=w_0EI/qL^4$).}
    \label{tab: eringen_validation}
    \begin{tabular}{c|c c}
    \hline\hline
    \multirow{2}{2em}{$h/l_f$} & \multicolumn{2}{c}{$\overline{w}$}\\
    \cline{2-3}
    & Present &  \cite{norouzzadeh2017finite,fernandez2016bending} \\
    \hline
       5  & 0.1169 & 0.1179\\
       10 & 0.1131 & 0.1139\\
       15& 0.1119 &0.1122\\
       20& 0.1113 & 0.1115\\
       \hline\hline
    \end{tabular}
\end{table}

\subsection{Convergence}
\label{ssec:convergence}
\begin{figure}[ht!]
    \begin{subfigure}[t]{0.5\textwidth}
    \centering
    \includegraphics[width=\textwidth]{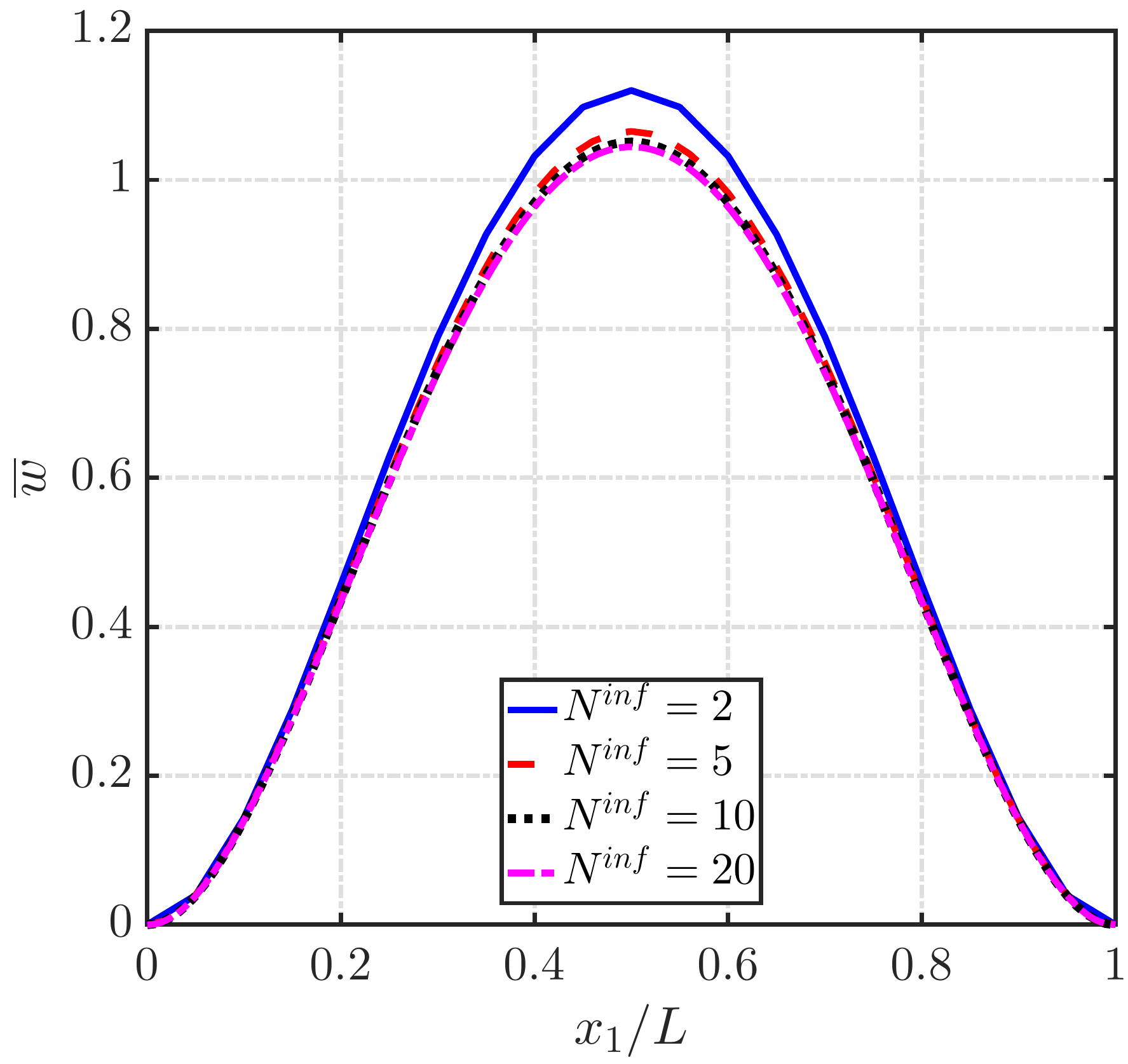}
        \caption{Two-noded $C^1$ element.}
        \label{fig: h-convergence}
    \end{subfigure}%
    ~ 
    \begin{subfigure}[t]{0.5\textwidth}
        \centering
    \includegraphics[width=\textwidth]{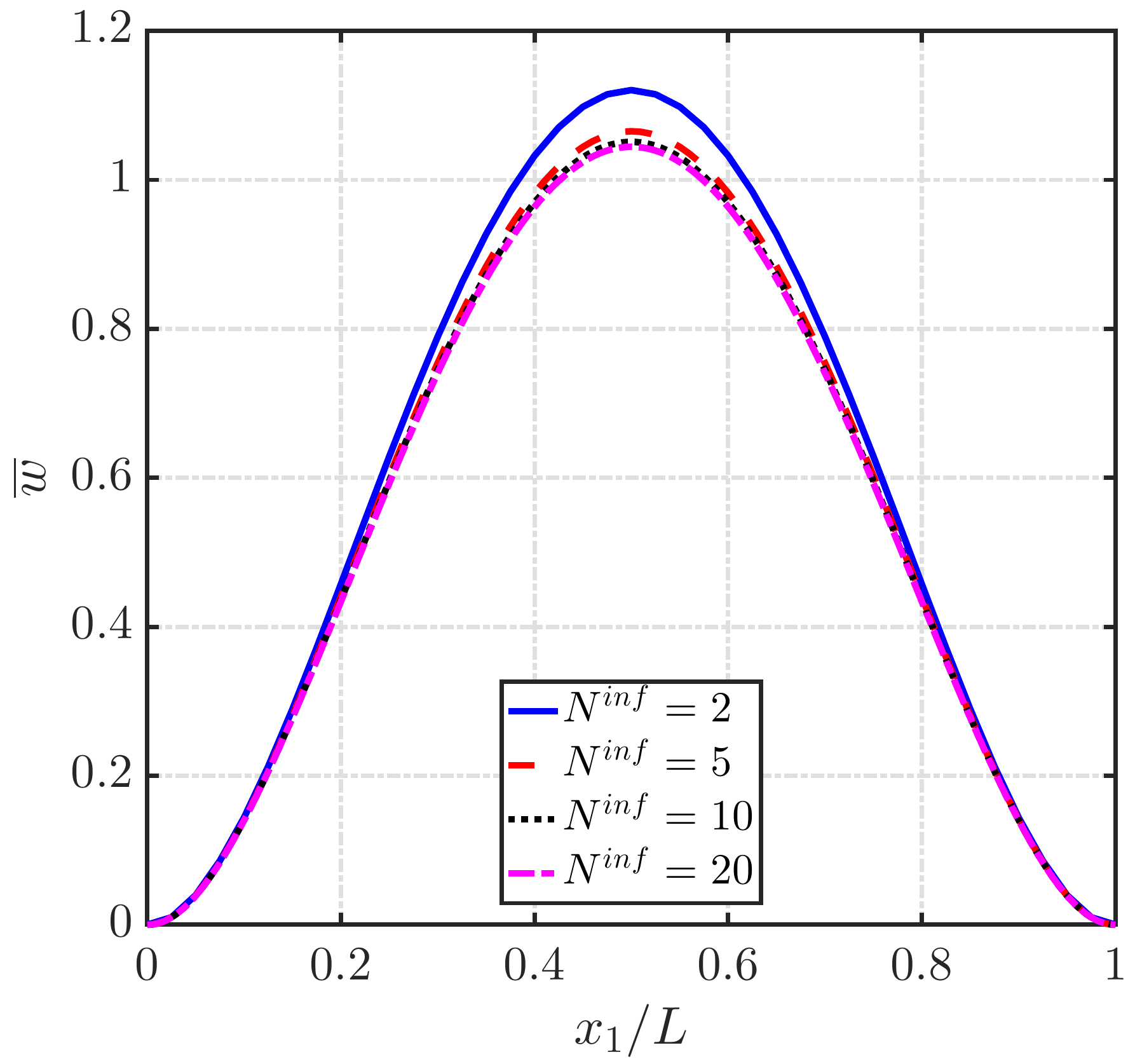}
        \caption{Three-noded $C^1$ element.}
        \label{fig: p-convergence}
    \end{subfigure}
    \caption{Numerical convergence of the f-FEM for $h$- and $p-$refinement ($\alpha=0.8$, $l_f=L/10$).}
    \label{fig: convergence}
\end{figure}

We establish the convergence of the f-FEM formulation by performing both $h-$ and $p-$ refinements. Any discretization in the mesh would result in larger number of elements being included in the horizon of nonlocality of the fractional-order model. The increased resolution would result in lower inconsistencies due to the truncation of the nonlocal horizon caused due to the ceil and floor operations described in \S\ref{numerical_integration_scheme}. Therefore, convergence is expected for the case where sufficient number of elements have been included in the zones of nonlocal interaction. This is established from the numerical results for the f-FEM of the beam in Fig.~(\ref{fig: convergence}). The convergence of the results obtained for increasing $N^{inf}$ ($=l_f/l_e$), referred to as the “dynamic rate of convergence” \cite{norouzzadeh2017finite}, indicate the necessary level of discretization for the convergence of f-FEM solution. 
Additionally, a three-noded $C^1$ model has been used to study the improvements obtained following the use of higher-order polynomials as shape functions. The numerical results in Fig.~(\ref{fig: p-convergence}) indicate convergence for this choice of the element. The convergence for $h-$ and $p-$ refinements indicate the consistency of the above developed f-FEM for modeling fractional-order nonlocal elastic solids.

A comprehensive evaluation of the convergence of the f-FEM is presented in Table \ref{tab: convergence_comb}. Both two- and three-noded elements as well as different combinations of the fractional parameters $\alpha$ and $l_f$ are considered. More specifically, the normalized maximum transverse displacement of the beam ($\overline{w})$ is compared for a given fractional order and length scale (compare the results by moving from top to bottom within a column in Table \ref{tab: convergence_comb}).
Values show that by increasing $l_f$ and reducing $\alpha$, that is by increasing the degree of nonlocality, larger number of elements are necessary for the convergence of f-FEM. 
It appears from Table \ref{tab: convergence_comb} that targeting an error threshold less than 1\% between successive refinements, the dynamic rate of convergence is $N^{inf}=10$. Therefore, the f-FEM simulations presented in the manuscript were performed at this level of discretization.

\begin{table}[h!]
    \centering
    \begin{tabular}{c | c |c c c c | c c c c}
    \hline\hline
       \multirow{3}{3em}{~~~$l_f$} & \multirow{3}{3em}{~$N^{inf}$} & \multicolumn{8}{c}{$\overline{w}$}\\
       \cline{3-10}
       &  & \multicolumn{4}{c|}{Two-noded}& \multicolumn{4}{c}{Three-noded}\\
         \cline{3-10}
         &  & $\alpha=1.0$& $\alpha=0.9$& $\alpha=0.8$ & $\alpha=0.7$ & $\alpha=1.0$& $\alpha=0.9$& $\alpha=0.8$ & $\alpha=0.7$ \\
         \hline\hline
          \multirow{4}{5em}{$L/5$} & 2 & 1.0000 & 1.0878 & 1.1776 & 1.2803 & 1.0000 & 1.0868 & 1.1773 & 1.2862\\
        & 5 & 1.0000 & 1.0750 & 1.1445 & 1.2211 & 1.0000 & 1.0743 & 1.1468 & 1.2251\\
          & 10 & 1.0000 & 1.0730 & 1.1424 & 1.2140 & 1.0000 & 1.0719 & 1.1410 & 1.2132\\
          & 20 & 1.0000 & 1.0720 & 1.1401 & 1.2098 & 1.0000 & 1.0709 & 1.1388 & 1.2088\\
          \hline
         \multirow{4}{5em}{$L/10$} & 2 & 1.0000 & 1.0602 & 1.1218 & 1.1907 & 1.0000 & 1.0590 & 1.1205 & 1.1914\\
         & 5 & 1.0000 & 1.0344 & 1.0667 & 1.1007 & 1.0000 & 1.0332 & 1.0653 & 1.0999\\
         & 10 & 1.0000 & 1.0275 & 1.0523 & 1.0778 & 1.0000 & 1.0264 & 1.0510 & 1.0768\\
         & 20 & 1.0000 & 1.0243 & 1.0456 & 1.0673 & 1.0000 & 1.0232 & 1.0444 & 1.0663\\
         \hline
          \multirow{4}{5em}{$L/20$} & 2 & 1.0000 & 1.0577  & 1.1172 & 1.1815 & 1.0000 & 1.0565  & 1.1155 & 1.1805\\
          & 5 & 1.0000 & 1.0244 & 1.0476 & 1.0717 & 1.0000 & 1.0232 & 1.0461 & 1.0705\\
          & 10 & 1.0000 & 1.0153 & 1.0288 & 1.0429 & 1.0000 & 1.0141 & 1.0275 & 1.0418\\
          & 20 & 1.0000 & 1.0109 & 1.0200 & 1.0294 & 1.0000 & 1.0098 & 1.0188 & 1.0285\\
          \hline\hline
    \end{tabular}
    \caption{Convergence of the f-FEM for a clamped-clamped beam for different fractional parameters.}
    \label{tab: convergence_comb}
\end{table}

\subsection{Effect of fractional constitutive parameters}
\label{ssec:effect_of_fractional_parameters}
\begin{figure*}[ht!]
    \centering
    \begin{subfigure}[t]{0.5\textwidth}
        \centering
        \includegraphics[width=\textwidth]{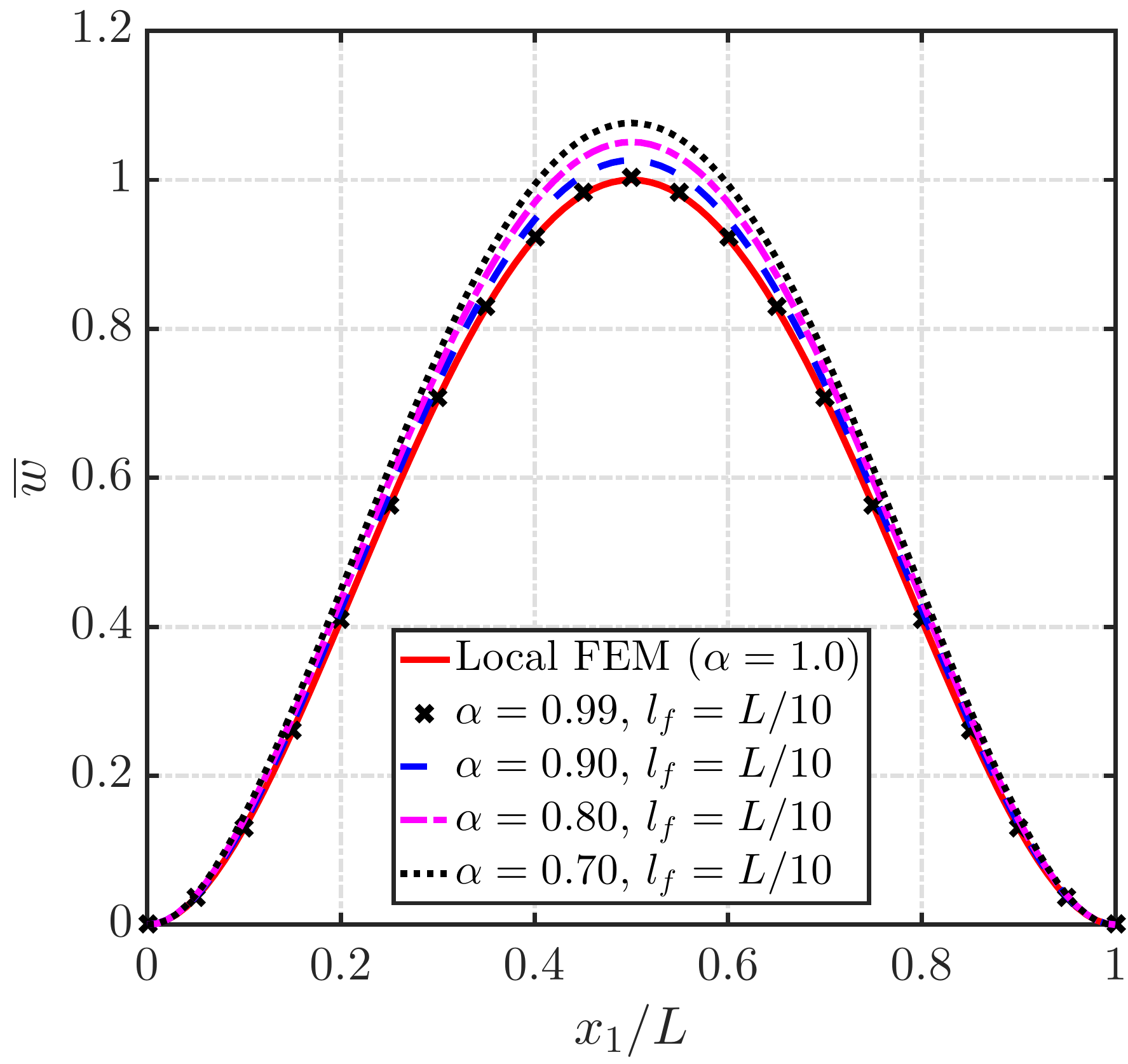}
        \caption{$\overline{w}$ vs $\alpha$.}
        \label{fig: cc_w_alpha}
    \end{subfigure}%
    ~ 
    \begin{subfigure}[t]{0.5\textwidth}
        \centering
        \includegraphics[width=\textwidth]{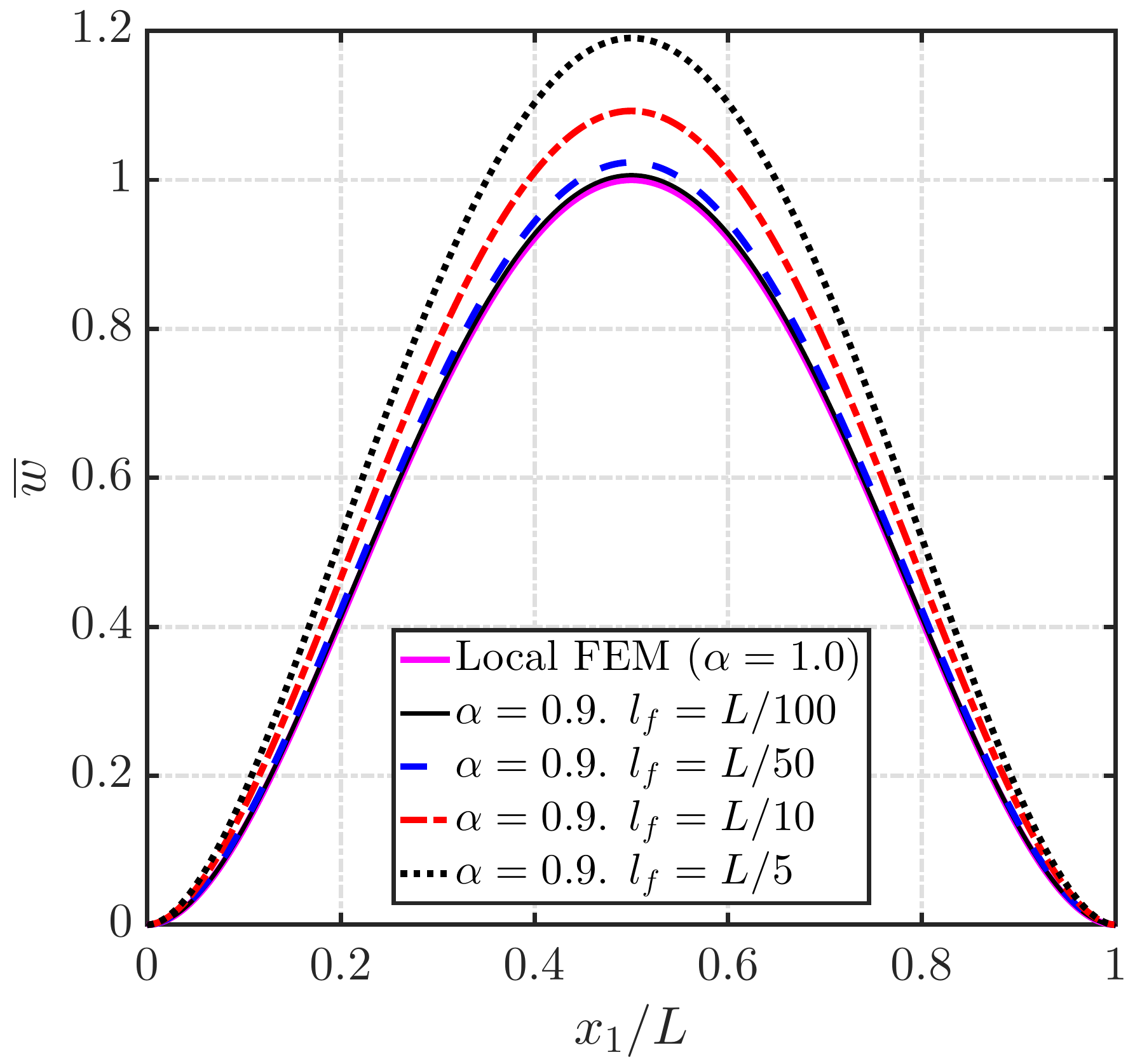}
        \caption{$\overline{w}$ vs $l_f$.}
        \label{fig: cc_w_lf}
    \end{subfigure}
    \caption{Effect of the fractional-order constitutive properties over the nonlocal elastic transverse displacement of the clamped beam.}
    \label{fig: cc_w}
\end{figure*}

\begin{figure*}[ht!]
    \centering
    \begin{subfigure}[t]{0.5\textwidth}
        \centering
        \includegraphics[width=\textwidth]{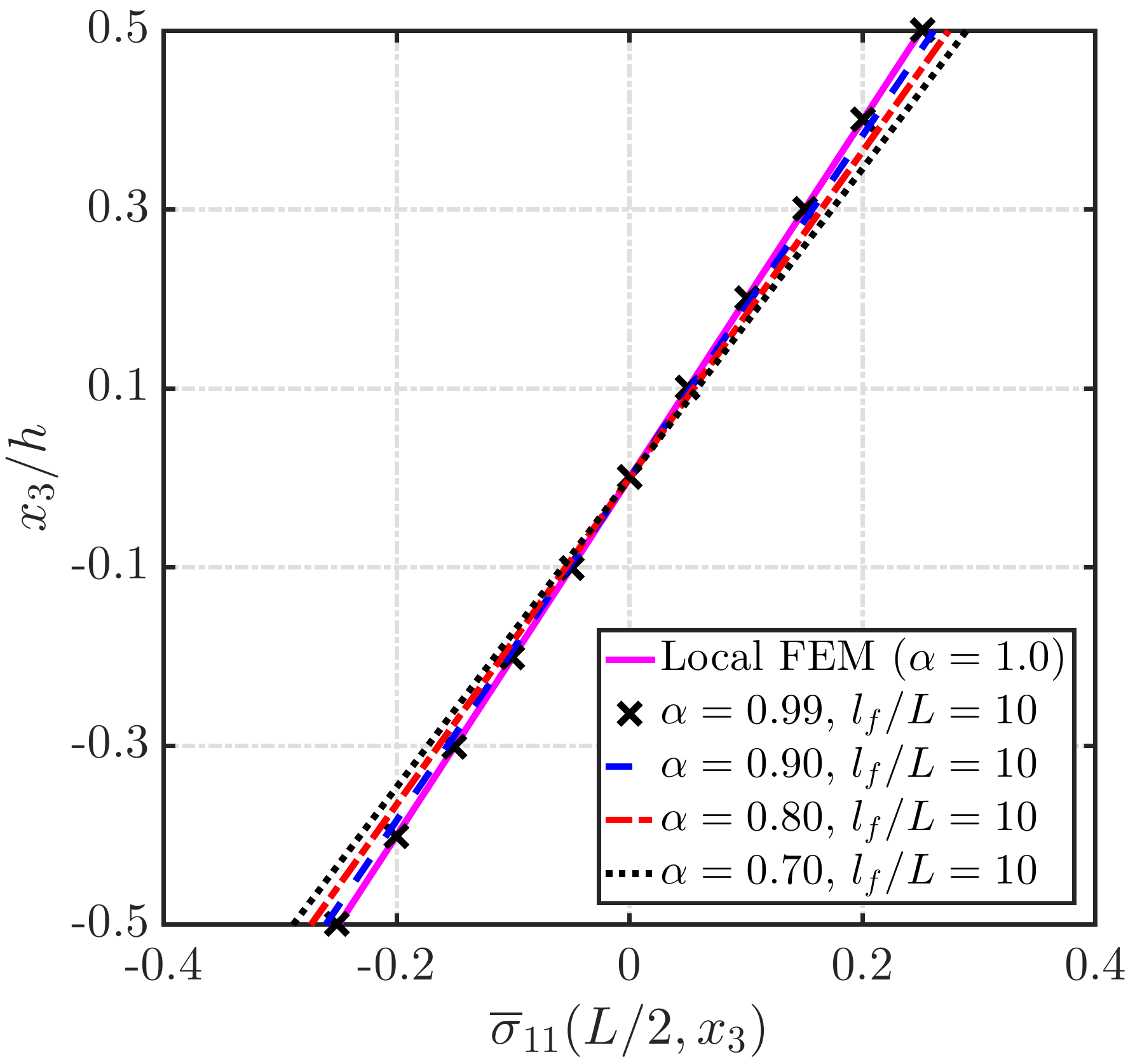}
        \caption{$\overline{\sigma}_{11}$ vs $\alpha$.}
        \label{fig: cc_sig11_alpha}
    \end{subfigure}%
    ~ 
    \begin{subfigure}[t]{0.5\textwidth}
        \centering
        \includegraphics[width=\textwidth]{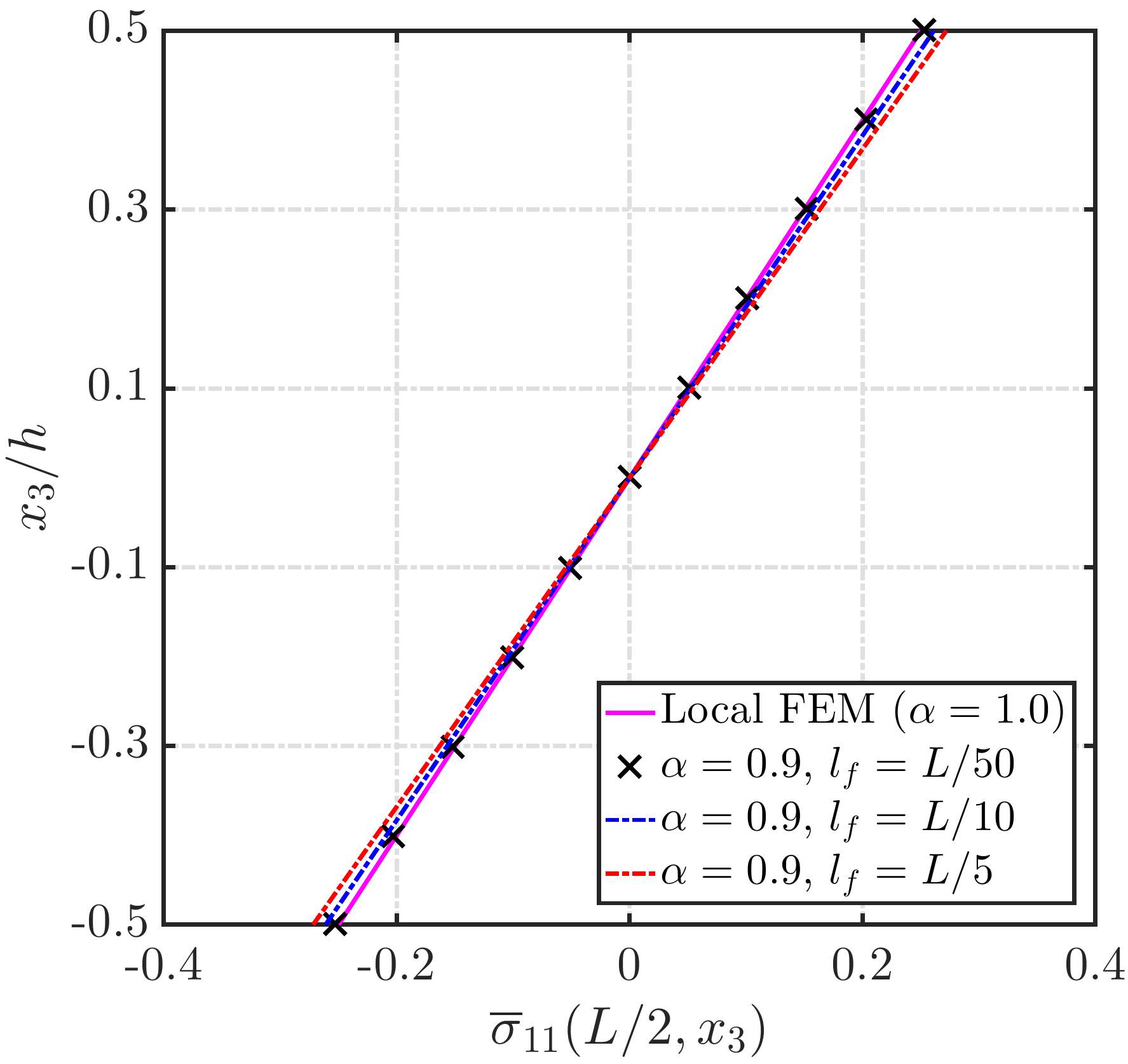}
        \caption{$\overline{\sigma}_{11}$ vs $l_f$.}
        \label{fig: c_sig11_lf}
    \end{subfigure}
    \caption{Effect of the fractional-order constitutive properties over the nonlocal elastic axial normal stress of the clamped beam.}
    \label{fig: cc_sig11}
\end{figure*}

Having validated the f-FEM and established its consistency, we now use the f-FEM to analyze the static response of the fractional-order nonlocal Euler-Bernoulli beam. More specifically, we analyze the effect of the fractional model parameters $\alpha$ and $l_f$ on the response of the fractional-order beam under different loading conditions. We start by considering a beam that is clamped at both ends and is subject to a UDL of magnitude $q_0$ (in N/m). The transverse displacement and the axial stress at the mid-plane of the beam are obtained for different values of $\alpha$ and $l_f$. The transverse displacement and the axial stress for the clamped beam are normalized in the following manner \cite{khodabakhshi2015unified}:
\begin{subequations}
\label{eq: norm_cc}
\begin{equation}
\overline{w}(x_1)=\frac{384~E~I}{q_0~L^4}w_0(x_1)
\end{equation}
\begin{equation}
\overline{\sigma}_{11}(L/2,x_3)=\frac{1}{q_0}\left(\frac{h}{L}\right)^2\tilde{\sigma}_{11}(L/2,x_3)
\end{equation}
\end{subequations}
The results of this study are shown in Fig.~(\ref{fig: cc_w}). As evident from Fig.~(\ref{fig: cc_w_alpha}), an increase in the transverse displacement is noted for reducing values of $\alpha$. This reduction in the stiffness upon considering the fractional-order nonlocal elasticity is in agreement with previous observations pertaining to nonlocal elasticity studies in literature. Similarly, an increase in the transverse displacement is observed from Fig.~(\ref{fig: cc_w_lf}) with the increasing horizon of nonlocality. The decrease in the stiffness of the structure due to the increasing degree of nonlocality (by reducing $\alpha$ and/or increasing $l_f$) is also established from the increase in the axial normal stress $\bar{\sigma}_{11}$ across the thickness at mid-length of the beam (see Fig. \ref{fig: cc_sig11}). The higher values of $\bar{\sigma}_{11}$ indicate the higher bending caused by reduced stiffness due to the fractional-order nonlocal elastic behaviour, resulting in higher values for the transverse displacement $\bar{w}$ when compared against the response of a local beam (see Fig.~(\ref{fig: cc_w})).

\begin{figure*}[t!]
    \centering
    \begin{subfigure}[t]{0.5\textwidth}
        \centering
        \includegraphics[width=\textwidth]{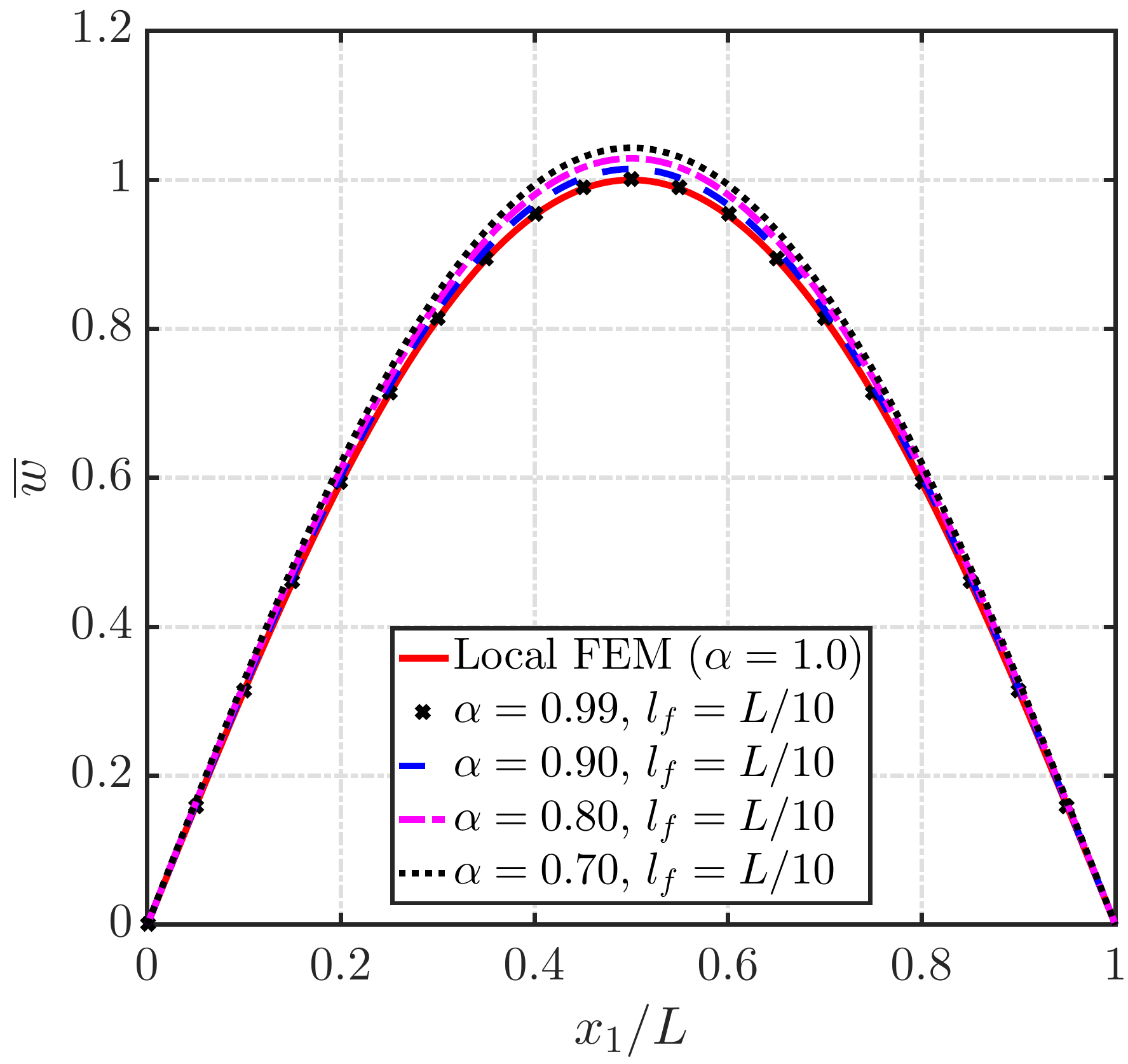}
        \caption{$\overline{w}$ vs $\alpha$.}
        \label{fig: ss_w_alpha}
    \end{subfigure}%
    ~ 
    \begin{subfigure}[t]{0.5\textwidth}
        \centering
        \includegraphics[width=\textwidth]{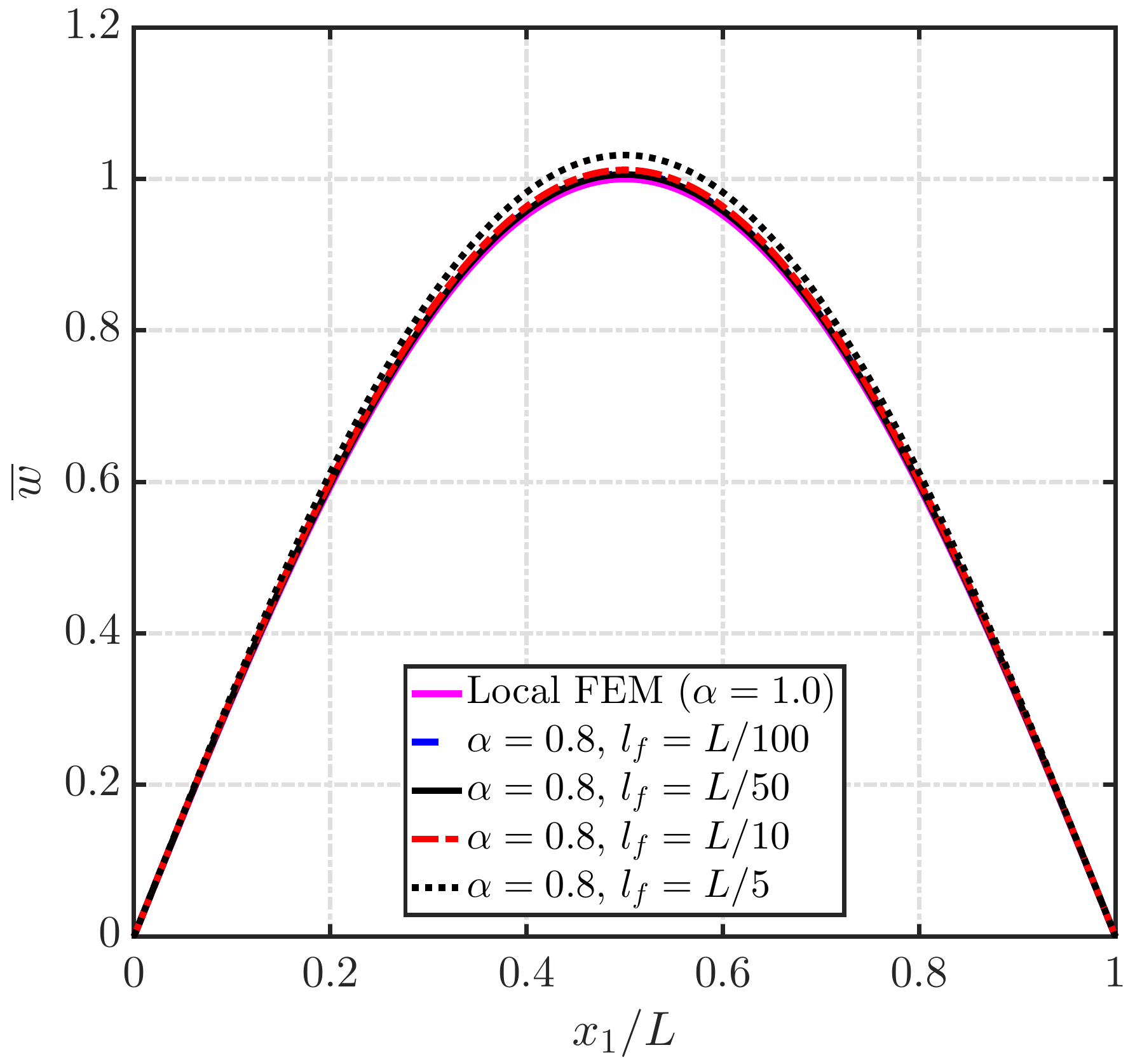}
        \caption{$\overline{w}$ vs $l_f$.}
        \label{fig: ss_w_lf}
    \end{subfigure}
    \caption{Effect of the fractional-order constitutive properties over the nonlocal elastic response of the simply-supported beam.}
    \label{fig: ss_w}
\end{figure*}

\begin{figure*}[t!]
    \centering
    \begin{subfigure}[t]{0.5\textwidth}
        \centering
        \includegraphics[width=\textwidth]{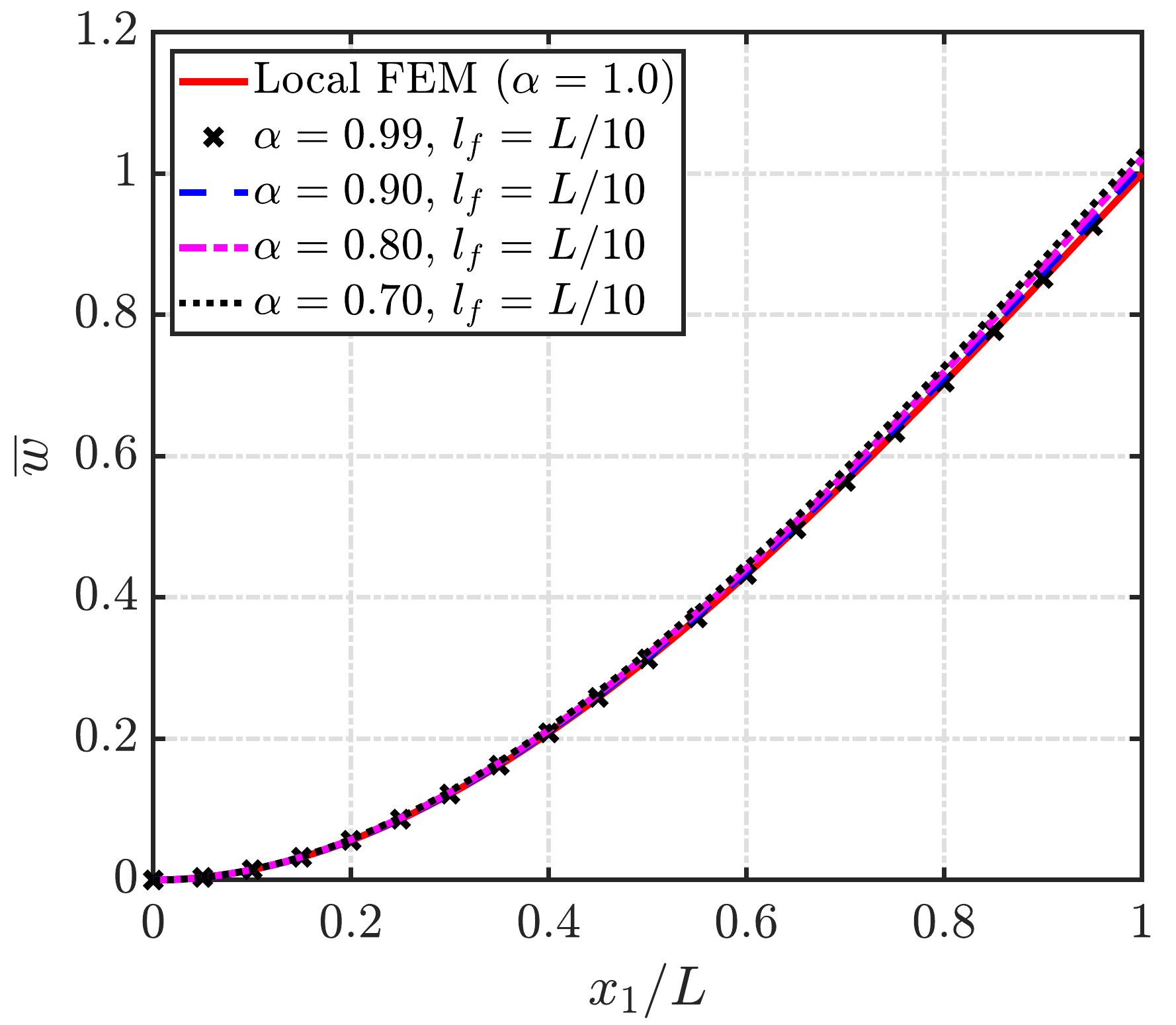}
        \caption{$\overline{w}$ vs $\alpha$.}
        \label{fig: cf_w_alpha}
    \end{subfigure}%
    ~ 
    \begin{subfigure}[t]{0.5\textwidth}
        \centering
        \includegraphics[width=\textwidth]{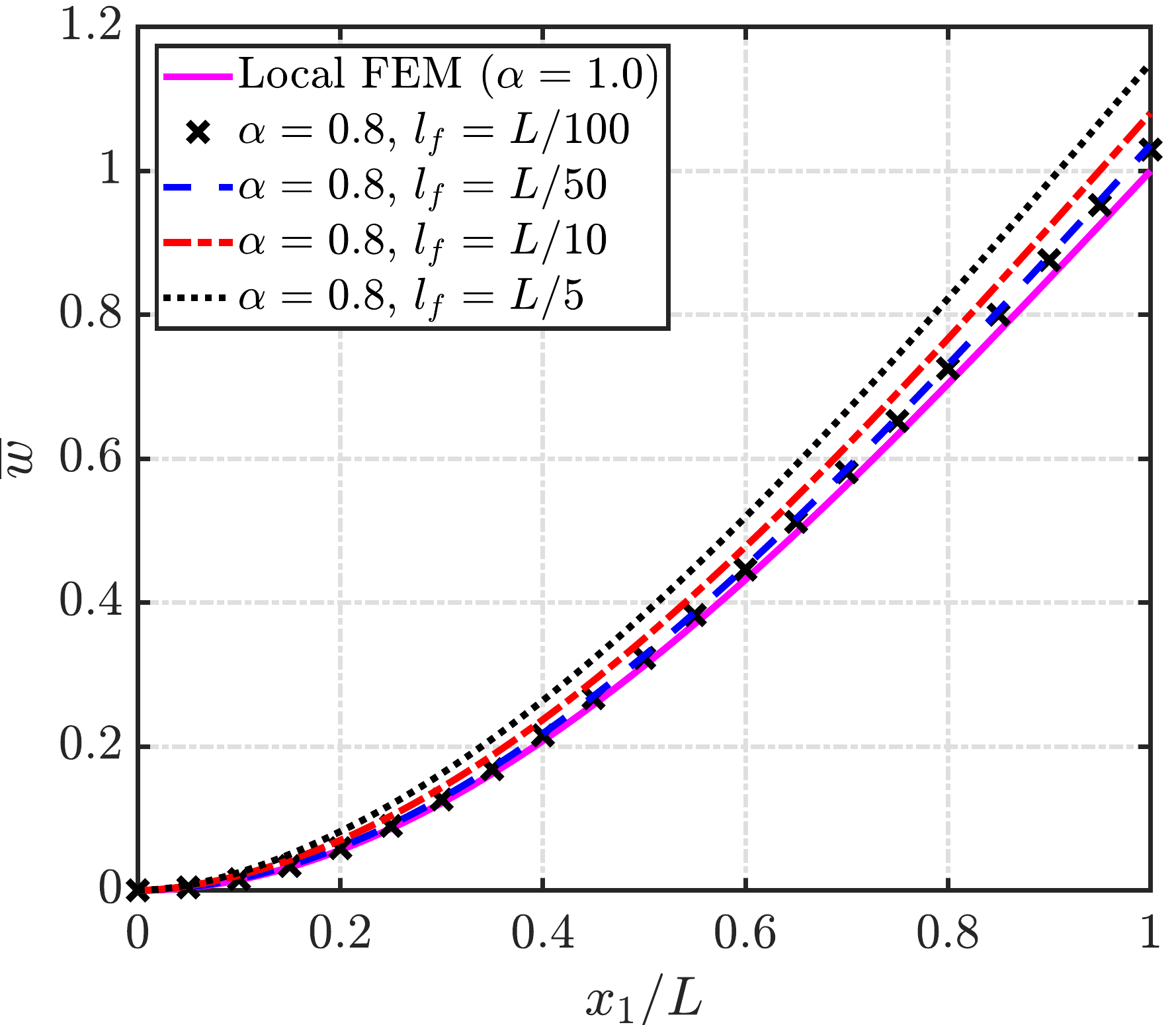}
        \caption{$\overline{w}$ vs $l_f$.}
        \label{fig: cf_w_lf}
    \end{subfigure}
    \caption{Effect of the fractional-order constitutive properties over the nonlocal elastic response of the cantilever beam.}
    \label{fig: cf_w}
\end{figure*}

We also studied the response of a simply-supported (S-S) beam and a cantilever (C-F) beam for varying fractional model parameters $\alpha$ and $l_f$. The normalized transverse displacements along the length of the beam for these boundary conditions are presented in Figs.~(\ref{fig: ss_w},\ref{fig: cf_w}). For these boundary conditions, the transverse displacement across the length has been normalized as follows:
\begin{subequations}
\begin{equation}
    \text{S-S:~~}\overline{w}(x_1)=\frac{384EI}{5q_0L^4}{w_0}(x_1)
\end{equation}
\begin{equation}
    \text{C-F:~~}\overline{w}(x_1)=\frac{3EI}{PL^3}{w_0}(x_1)
\end{equation}
\end{subequations}
where, $q_0$ is the magnitude of the UDL applied over S-S beam, and $P$ is the magnitude of the concentrated load applied at the tip of the cantilever beam. The observations noted previously for the clamped beam are also noted for the beam with the above boundary conditions. Thus the nonlocal interactions reduce the stiffness of the structure irrespective of the boundary conditions. 

The consistency of these results, obtained from the fractional-order modeling for the strain-displacement relations, indicates the consistency and complete nature of modeling nonlocal elastic interactions by the fractional-order derivatives. In contrast to classical approaches to nonlocal elasticity, no paradoxical results are obtained for different loading and boundary conditions from the fractional-order model of nonlocal elasticity. More specifically, this is in contrast to the prediction of hardening or an absence of nonlocal interactions altogether, for cantilever beams subjected to a point load at the free end, as noted from the differential model of Eringen's nonlocal elasticity \cite{challamel2008small}. Similar paradoxes have also been noted for the case of a simply-supported nonlocal beam subjected to a UDL upon using a two-phase integro-differential model of nonlocal elasticity \cite{khodabakhshi2015unified}. 
Further, the general nature of the fractional-order model for nonlocal elasticity developed here becomes clear from a comparison of various nonlocal theories like the two-phase nonlocal integral constitutive model \cite{polizzotto2001nonlocal}, modifications in the kernel function corresponding to the convolution of nonlocal interactions \cite{koutsoumaris2017different}, and combined nonlocal and gradient elasticity theories \cite{benvenuti2013one}. In contrast to the above mentioned studies, the definition of the nonlocal modeling based on the RC fractional derivative employed here, is uniform for all the cases. Moreover, the necessity and flexibility of the f-FEM developed here, in a study of integral model of nonlocal elasticity, becomes clear from the complexities involved in the handling of integral boundary conditions \cite{fernandez2016bending}.


We make an important remark concerning the physically acceptable range for the order $\alpha$. As demonstrated above, the degree of nonlocality increases  with decreasing $\alpha$ leading to a consistent softening of the structure. However, as shown in \cite{sumelka2014thermoelasticity,sumelka2015non}, results for very low values of $\alpha$ ($\approx0.2$, which indicates a very strong nonlocality in the fractional sense) lead to non-physical solutions. Hence, there exists a limit on the order of the RC fractional derivative \cite{sumelka2015non}. In other terms, it can be concluded that the fractional-calculus based homogenization or modeling of nonlocality breaks down for values of alpha close to the lower integer limit. While, in the present model, the specific value is $\alpha \approx 0.2$, we note that this threshold value can vary as a result of different choices of the model's parameters (particularly, the size of the nonlocal horizon). We emphasize that this breakdown is not a characteristic of the f-FEM technique as the same observation is also noted when using finite difference methods to obtain the numerical solutions (see, for example, \cite{sumelka2014thermoelasticity,sumelka2015non}).

We present a possible interpretation to the breakdown of the fractional-order nonlocal model beyond a certain order close to the lower integer limit. Note that a fractional-order derivative can be considered as an interpolation between consecutive integer-order derivatives ($n-1<\alpha<n$; $n\in\mathrm{N}^+$). It follows that also the physical mechanisms simulated via such derivatives are combined into a hybrid mechanism \cite{podlubny1998fractional}. A typical example consists in the viscoelastic behavior which is typically modeled using a time-fractional derivative of displacement with $\alpha=0.5$. Such derivative provides an interpolation between the elastic (modeled using a time-fractional derivative of displacement with $\alpha=0$) and the viscous (modeled using a time-fractional derivative of displacement with $\alpha=1$) behaviors. 
Further, the contributions of the integer-order derivatives (calculated at the bounding integer values of $\alpha$) to the final value of the fractional-order derivative depend on the specific values of both the order and the interval of the derivative. In other terms, for a fixed interval and for $n-1<\alpha<n$ ($n\in\mathrm{N}^+$), the contribution of the $n-1^{th}$ integer-order derivative would be higher if $\alpha$ is closer to $n-1$ instead of $n$.

In the context of the above discussion, note that the strain in the fractional-order formulation can be considered as an interpolation between the displacement (zero-order derivative) and the gradient (first-order derivative) of displacement. From the constitutive relation in Eq.~(\ref{stress_equation}), it follows that the stress at a point (and consequently, the force experienced by the point) is directly proportional to both the displacement and the gradient (first-order) of displacement, via a positive constant of proportionality. Recall that when the force of a point is directly proportional to the displacement of a point (via a positive constant), the resulting motion is generally exponential and unbounded \footnote{$F\propto C^2u \implies \frac{d^2 u}{d t^2}=C^\dagger u \implies u = A_1 e^{C^\dagger t} + A_1 e^{-C^\dagger t}$}. On the contrary, when the force (or stress) is proportional to the displacement gradient the resulting motion is uniform and bounded (as in classical local continuum mechanics). It follows that there exists a critical order $\alpha$, beyond which the fractional-order interpolation of the displacement and the gradient (first-order) of displacement becomes unstable due to larger contribution of the displacement term. Note that the specific value of this critical $\alpha$ would depend on the specific interval-length chosen in the formulation. In the current study, for a maximum length-scale of $l_f = L/5$, this critical value of $\alpha\approx0.2$ which agrees with previously conducted studies \cite{sumelka2014thermoelasticity,sumelka2015non}. It was also observed that with increasing value of the length-scale $l_f$, the value of the critical order $\alpha$ gets closer to zero.

\section{Conclusions}

This paper presented a finite element model for the numerical solution of fractional-order boundary value problems. A fractional-order nonlocal Euler-Bernoulli beam based on a frame-invariant and dimensionally consistent fractional-order nonlocal continuum theory was considered. The governing equation of the beam under Euler-Bernoulli type constitutive modeling was derived in strong form using variational principles. We showed that the fractional-order nonlocal modeling results in a self-adjoint and positive definite system with a unique solution. The numerical model for this continuum system was developed using the minimum potential energy principle. For the fractional-order system, energy minimization was carried out on a global scale owing to the nonlocal nature of the response, hence resulting in pre-assembled system matrices for the f-FEM. Additionally, we presented a scheme to circumvent the singularity associated with the convolution kernel of the fractional-order model for nonlocal elasticity. This also represented a critical step in order to develop Ritz finite element formulations.
The proposed f-FEM model was validated with benchmark problems in fractional-order equations and nonlocal elasticity. Then, the nonlocal elastic response of a Euler-Bernoulli subjected to various loading and boundary conditions was considered. Paradoxical observations noted in literature, such as absence of nonlocal interactions and hardening for cantilever and simply-supported beams have been addressed with the use of the positive-definite fractional-order model for nonlocal elasticity. It is important to highlight that the f-FEM methodology developed is very general and, although in the context of this study it was tested for 1-D elasticity, it can be easily extended to plates and shells.

\section*{Appendix 1}
{\textbf{Proof of Frame Invariance:}}
Recalling the definition of the deformation gradient tensor, we must show that the individual terms in Eq.~(\ref{Fractional_F_net}) are frame invariant. We start by investigating the frame invariance of $\tilde{\textbf{F}}_{X}$.
Consider a rigid-body motion superimposed on a general point $\textbf{X}$ (see Fig.~(1b) of the manuscript) of the reference configuration of the body as:
\begin{equation}
\label{sq5}
\bm{\Psi}(\textbf{X},t)=\textbf{c}(t)+\textbf{Q}(t)\textbf{X},
\end{equation}
where $\textbf{Q}(t)$ is a proper orthogonal tensor denoting a rotation and $\textbf{c}(t)$ is a spatially constant term representing a translation. Under this rigid-body motion, the fractional deformation gradient $\tilde{\textbf{F}}^{\Psi}_X$ should be an orthogonal tensor such that $\tilde{\textbf{F}}^{\Psi T}_X\tilde{\textbf{F}}^{\Psi}_X = \textbf{I}$. More specifically, the fractional deformation gradient tensor should transform as $\tilde{\textbf{F}}^{\Psi}_X = \textbf{Q}$ (similar to the classical continuum case where $\textbf{F}^{\Psi} = \textbf{Q}$) such that the strain measures are null. From the definition of $\tilde{\textbf{F}}_X$  given in Eq.~(\ref{Fractional_F}a) of our study it follows that:
\begin{equation}
\label{sq6}
\begin{split}
\tilde{\textbf{F}}_{X_{ij}}^{\Psi} = \frac{1}{2} \Gamma(2-\alpha) \biggl[ \frac{L_{A_{j}}^{\alpha-1}}{\Gamma(1-\alpha)} \int_{X_{A_{j}}}^{X_j} \frac{D^1_{S_j} \Psi_{i}(\textbf{S},t)}{(X_j-S_j)^{\alpha}}  \mathrm{d}S_j + \frac{L_{B_{j}}^{\alpha-1}}{\Gamma(1-\alpha)} \int_{X_{j}}^{X_{B_j}} \frac{D^1_{S_j} \Psi_{i}(\textbf{S},t)}{(S_j-X_j)^{\alpha}} \mathrm{d}S_j \biggr]
\end{split}
\end{equation}
where $\textbf{S}$ is a dummy vector representing the spatial variable. Further, $D^1_{S_j} \Psi_{i}(\textbf{S},t)$ simplifies as:
\begin{equation}
\label{sq7}
D^1_{S_j} \Psi_{i}(\textbf{S},t)=\frac{\mathrm{d}\Psi_{i}(\textbf{S},t)}{\mathrm{d}S_j}=\frac{\mathrm{d}}{\mathrm{d}S_j}(c_i+Q_{ik}S_{k})
\end{equation}
Noting that $\frac{\mathrm{d}c_i(t)}{\mathrm{d}S_j}=0$ and $\textbf{Q}=\textbf{Q}(t)$ it follows that:
\begin{equation}
\label{sq8}
D^1_{S_j} \Psi_{i}(\textbf{S},t)=Q_{ik}S_{k,j}=Q_{ik}\delta_{kj}=Q_{ij}
\end{equation}
Thus, under the rigid body motion $\bm{\Psi}$:
\begin{equation}
\label{sq9}
\tilde{\textbf{F}}_{X_{ij}}^{\Psi}=\frac{1}{2}\bigl[L_{A_{j}}^{\alpha-1} (X_j-X_{A_{j}})^{1-\alpha}+L_{B_{j}}^{\alpha-1} (X_{B_{j}}-X_j)^{1-\alpha}\bigr]Q_{ij}
\end{equation}
In the above simplifications we have used the following property of the $\Gamma(\cdot)$ function: $\Gamma(2-\alpha) = (1-\alpha)\Gamma(1-\alpha)$. As highlighted in the study, the length-scales $L_{A_j}$ and $L_{B_j}$ are taken such that: $L_{A_{j}}=X_j-X_{A_{j}}$ and $L_{B_{j}}=X_{B_{j}}-X_j$. This has also been illustrated schematically in Fig.~(\ref{fig1}b). By substituting these relations in Eq.~(\ref{sq9}), it follows that $\tilde{\textbf{F}}^{\Psi}_X = \textbf{Q}$ at all times. We also emphasize that the nonlocal formulation allows for an exact treatment of frame invariance in the presence of asymmetric horizons which occurs at points close to material boundaries and interfaces. The different horizon lengths $L_{A_j}$ and $L_{B_j}$ enables the truncation of the horizon at points close to or on the boundary in order to exactly satisfy frame-invariance. By repeating this procedure, it is immediate to show that the same arguments hold for the frame invariance of $\tilde{\textbf{F}}^{\Psi}_x$ and subsequently of $\overset{\alpha}{{\textbf{F}}}$.

\section*{Appendix 2}
{\textbf{Relation of Fractional-Order Continuum Model to Eringen's Model:}} In the following, we demonstrate that a specialization of Eringen's integral approach to nonlocality by means of a power-law attenuation kernel leads to fractional-order constitutive relations with Caputo derivatives. More specifically, our approach makes use of a Riesz-Caputo definition of the fractional operator that uses proper length scales to guarantee both dimensional consistency and frame invariance.

For the sake of this discussion, consider a 1D bar of length $L$. Eringen's integral stress-strain constitutive relation at a point $x$ inside the bar ($x\in[0,L]$) is given by:
\begin{equation}
\label{eq: constitutive_eringen_step1}
    \sigma(x) = E\left[ \int_a^b \mathcal{K}(x-s) \epsilon(s) \mathrm{d}s \right]
\end{equation}
where $\sigma(x)$ and $\epsilon(x)$ denote the stress and strain at a specific point $x$ in the bar and $E$ denotes the Young's modulus of the solid. The interval $[a,b]$ corresponds to the horizon of nonlocality. Note that, when the domain of the bar is finite, the horizon of nonlocality has to be within the finite domain of the structure, that is $[a,b]\subset[0,L]$. $\mathcal{K}$ is the kernel corresponding to the strength of the nonlocal interactions. For the 1D bar, the above stress-strain relation can be expressed as a function of the displacement field $u$ as:
\begin{equation}
\label{eq: constitutive_eringen_step2}
    \sigma(x) = E\left[ \int_a^b \mathcal{K}(x-s) \frac{\mathrm{d}u(s)}{\mathrm{d}s} \mathrm{d}s \right]
\end{equation}
Now, consider the following selection for the nonlocal kernel:
\begin{equation}
\label{eq: kernel}
    \mathcal{K}(x,s)=\left\{\begin{matrix}
    \frac{1}{2}(1-\alpha)l_A^{\alpha-1}{(x-s)^{-\alpha}} ~ \forall~ s\in{(a,x)}\\
    \frac{1}{2}(1-\alpha)l_B^{\alpha-1}{(s-x)^{-\alpha}} ~ \forall~ s\in{(x,b)}
    \end{matrix}\right.
\end{equation}
where $l_A$ and $l_B$ are length scales such that $l_A=x-a$ and $l_B=b-x$, and $\alpha$ is a non-negative parameter such that $\alpha\in(0,1)$. Note that the kernel is positive definite in nature which is a pre-requisite for a stable nonlocal formulation. Additionally, there are some important physical implications from the above definition for the nonlocal kernel. We have shown in our study that the above definition for the nonlocal kernel maintains dimensional consistency of the formulation. Also, the factors $\frac{1}{2}\Gamma(2-\alpha)$, $l_{A}^{\alpha-1}$, and $l_{B}^{\alpha-1}$ allow ensuring the frame invariance of the constitutive relations (see Appendix 1).

The generalization of the nonlocal kernel to Eq.~(\ref{eq: kernel}) leads to the following expression for the stress in the bar:
\begin{equation}
\label{eq: fractional_constitutive_step1}
    \sigma(x) = E\left[ \frac{1}{2}(1-\alpha)l_A^{\alpha-1} \int_a^x \frac{D^1_s u(s)}{(x-s)^\alpha} \mathrm{d}s + \frac{1}{2}(1-\alpha)l_B^{\alpha-1} \int_x^b \frac{D^1_s u(s)}{(s-x)^\alpha} \mathrm{d}s \right]
\end{equation}
where we have denoted the first integer-order derivative of the displacement field ${\mathrm{d}u(s)}/{\mathrm{d}s}$ as $D^1_s u(s)$. The above expression can be recast into the following form:
\begin{equation}
\label{eq: fractional_constitutive_step2}
    \sigma(x) = E\Bigg[ \underbrace{\frac{1}{2}\Gamma(2-\alpha) \Bigg\{ l_A^{\alpha-1} \underbrace{ \left[ \frac{1}{\Gamma(1-\alpha)} \int_a^x \frac{D^1_s u(s)}{(x-s)^\alpha} \mathrm{d}s\right]}_{\text{Left-handed Caputo derivative }} - l_B^{\alpha-1} \underbrace{\left[ - \frac{1}{\Gamma(1-\alpha)} \int_x^b \frac{D^1_s u(s)}{(s-x)^\alpha} \mathrm{d}s\right]}_{\text{Right-handed Caputo derivative }} \Bigg\}}_{\text{Riesz-Caputo derivative of the displacement}} \Bigg]
\end{equation}
where $\Gamma(\cdot)$ is the Gamma function. In the above mathematical manipulations, we have used the standard relation: $\Gamma(2-\alpha)=(1-\alpha)\Gamma(1-\alpha)$. From the above, we see that with our choice of the kernel the left- and right-handed Caputo derivatives appear in the nonlocal constitutive relations as a generalization of Eringen's integral approach. In our study, the combination of the left- and right-handed Caputo derivatives was referred to as the RC derivative (Eq.~(\ref{RC_definition})). By using the definition for the RC derivative, the constitutive relation in Eq.~(\ref{eq: fractional_constitutive_step2}) can now be expressed as:
\begin{equation}
\label{eq: fractional_constitutive_step3}
    \sigma(x) = E \underbrace{[D^\alpha_x u(x)]}_{{\substack{\text{Nonlocal}\\\text{strain}}}} = E\tilde{\epsilon}(x)
\end{equation}
The above equation is exactly the constitutive relation used in our study (Eq.~(\ref{eq: constt_axial})); in fact, $[D^\alpha_x u(x)]$ is the fractional-order strain $\tilde{\epsilon}(x)$. Note that we denoted the stress as $\tilde{\sigma}(x)$ in our study in order to be consistent with the use of a $\tilde{\square}$ for all physical quantities in the nonlocal solid. 
Note also that for $\alpha=1$, the classical constitutive relations for a 1D solid are obtained from Eq.~(\ref{eq: fractional_constitutive_step3}).

\section*{Appendix 3}
Consider the minimization integral $\delta\Pi$ in Eq.~(\ref{minimization_step_1}). The variations $\delta u_0(x_1)$ and $\delta w_0(x_1)$ in Eq.~(\ref{minimization_step_1}) are independent of each other. Hence the integrals in Eq.~(\ref{minimization_step_1}) are evaluated independently. For the sake of brevity, we outline only the steps involved in the mathematical operations over axial variables, i.e., Eq.~(\ref{minimization_step_2}a). These steps extend directly for the transverse displacement in Eq.~(\ref{minimization_step_2}b).
Using the definition of the Riesz-Caputo fractional derivative given in Eq.~(\ref{RC_definition}), the first integral corresponding to the variation of axial displacement in Eq.~(\ref{minimization_step_2}b) is obtained as:
\begin{equation}
\label{eq: var_axial}
\begin{split}
    \int_{0}^{L}N^\dagger(x_1)D_{x_1}^{\alpha}[\delta u_0(x_1)]\mathrm{d}x_1=\frac{1}{2}\Gamma(2-\alpha)\Bigg[l_A^{\alpha-1}\int_0^L N^\dagger(x_1)~{}^{C}_{x_1-l_A}D^{\alpha}_{x_1}[\delta u_0(x_1)]\mathrm{d}x_1-\\
    l_B^{\alpha-1}\int_0^L N^\dagger(x_1)~{}^{C}_{x_1}D^{\alpha}_{x_1+l_B}[\delta u_0(x_1)]\mathrm{d}x_1\Bigg]
\end{split}
\end{equation}
From the definitions for the left and right Caputo derivatives \cite{podlubny1998fractional} we obtain:
\begin{subequations}
\begin{equation}
\label{eq: left_caputo_axial}
     \int_0^L N^\dagger(x_1)~{}^{C}_{x_1-l_A}D^{\alpha}_{x_1}[\delta u_0(x_1)]\mathrm{d}x_1=\frac{1}{\Gamma(1-\alpha)}\int_0^L N^\dagger(x_1)\left[\int_{x_1-l_A}^{x_1}\left(x_1-s_1\right)^{-\alpha}\frac{\mathrm{d}\delta u_0(s_1)}{\mathrm{d}s_1}\mathrm{d}s_1\right]\mathrm{d}x_1
\end{equation}
\begin{equation}
\label{eq: right_caputo_axial}
    \int_0^L N^\dagger(x_1)~{}^{C}_{x_1}D^{\alpha}_{x_1+l_B}[\delta u_0(x_1)]\mathrm{d}x_1=-\frac{1}{\Gamma(1-\alpha)}\int_0^L N^\dagger(x_1)\left[\int_{x_1}^{x_1+l_B}\left(s_1-x_1\right)^{-\alpha}\frac{\mathrm{d}\delta u_0(s_1)}{\mathrm{d}s_1}\mathrm{d}s_1\right]\mathrm{d}x_1
\end{equation}
\end{subequations}
The above integrals are further evaluated using integration by parts in order to transfer the derivative from independent variable (displacement field) to the secondary variable (stress resultant). This leads to the following:
\begin{subequations}
\begin{equation}
    \int_0^L N^\dagger(x_1)~{}^{C}_{x_1-l_A}D^{\alpha}_{x_1}[\delta u_0(x_1)]\mathrm{d}x_1=\frac{1}{\Gamma(1-\alpha)}\int_0^L \frac{\mathrm{d}\delta u_0(s_1)}{\mathrm{d}s_1} \left[\int_{s_1}^{s_1+l_A}\left(x_1-s_1\right)^{-\alpha}N^\dagger(x_1)\mathrm{d}x_1\right]\mathrm{d}s_1
\end{equation}
\begin{equation}
   \int_0^L N^\dagger(x_1)~{}^{C}_{x_1}D^{\alpha}_{x_1+l_B}[\delta u_0(x_1)]\mathrm{d}x_1=-\frac{1}{\Gamma(1-\alpha)}\int_0^L \frac{\mathrm{d}\delta u_0(s_1)}{\mathrm{d}s_1} \left[\int_{s_1-l_B}^{s_1}\left(s_1-x_1\right)^{-\alpha}N^\dagger(x_1)\mathrm{d}x_1\right]\mathrm{d}s_1
\end{equation}
\end{subequations}
Using the definitions for left- and right- fractional integrals \cite{podlubny1998fractional} in the above results we obtain:
\begin{subequations}
\begin{equation}
\label{eq: left_caputo_axial_intpart1}
    \int_0^L N^\dagger(x_1)~{}^{C}_{x_1-l_A}D^{\alpha}_{x_1}[\delta u_0(x_1)]\mathrm{d}x_1=\int_0^L \left[\frac{\mathrm{d}\delta u_0(x_1)}{\mathrm{d}x_1}\right] {}_{x_1}I^{1-\alpha}_{x_1+l_A}[N^\dagger(x_1)]~\mathrm{d}x_1
\end{equation}
\begin{equation}
\label{eq: right_caputo_axial_intpart1}
   \int_0^L N^\dagger(x_1)~{}^{C}_{x_1}D^{\alpha}_{x_1+l_B}[\delta u_0(x_1)]\mathrm{d}x_1=-\int_0^L \left[\frac{\mathrm{d}\delta u_0(x_1)}{\mathrm{d}x_1}\right] {}_{x_1-l_B}I^{1-\alpha}_{x_1}[N^\dagger(x_1)]~\mathrm{d}x_1
\end{equation}
\end{subequations}
Repeating the integration by parts and substituting the resulting expressions in Eq.~(\ref{eq: var_axial}) we obtain:
\begin{equation}
\label{eq: var_axial_intermediate}
\begin{split}
    \int_{0}^{L}N^\dagger(x_1)D_{x_1}^{\alpha}[\delta u_0(x_1)]\mathrm{d}x_1= \frac{1}{2}\Gamma(2-\alpha)\bigg[l_A^{\alpha-1}\left.\left[{}_{x_1}I^{1-\alpha}_{x_1+l_A} N^\dagger(x_1)\delta u_0\right]\right\vert_{0}^{L}-
    l_B^{\alpha-1}\left.\left[{}_{x_1-l_B}I^{1-\alpha}_{x_1} N^\dagger(x_1)~\delta u_0\right]\right\vert_{0}^{L}\bigg]+\\
    \frac{1}{2}\Gamma(2-\alpha)\Bigg[l_A^{\alpha-1}\int_0^L\left[{}^{RL}_{x_1}D^{\alpha}_{x_1+l_A} N^\dagger(x_1)\right]\delta u_0(x_1)~\mathrm{d}x_1-l_B^{\alpha-1}\int_0^L\left[{}^{RL}_{x_1-l_B}D^{\alpha}_{x_1} N^\dagger(x_1)\right]\delta u_0(x_1)~\mathrm{d}x_1\Bigg]
\end{split}
\end{equation}
where ${}^{RL}_{x_1-l_B}D_{x_1}^\alpha(\cdot)$ and ${}^{RL}_{x_1}D_{x_1+l_A}^\alpha(\cdot)$ are the left- and right-handed Riemann Liouville derivatives, respectively. Note that in the fractional integral term $\left(l_A^{\alpha-1}{}_{x_1}I^{1-\alpha}_{x_1+l_A} N^\dagger(x_1)\right)$, in Eq.~(\ref{eq: var_axial_intermediate}) above, $l_A=0$ because $x_1\in\{0,l\}$, i.e., $x_1$ lies on either boundary of the beam. Similarly, $l_B=0$ in $\left(l_B^{\alpha-1}{}_{x_1-l_B}I^{1-\alpha}_{x_1} N^\dagger(x_1)\right)$. Under these limiting conditions the fractional integrals will converge to the function itself (similar to what is presented in Eqs.~(\ref{sq10}-\ref{sq13})) . Further, using the definition of Riesz RL derivative given in Eq.~\eqref{eq: riesz_rl_der} the above integral simplifies to Eq.~(\ref{minimization_step_2}b).
\section{Acknowledgements}
The following work was supported by the National Science Foundation (NSF) under the grants MOMS \#1761423 and DCSD \#1825837, and the Defense Advanced Research Project Agency (DARPA) under the grant \#D19AP00052. The content and information presented in this manuscript do not necessarily reflect the position or the policy of the government. The material is approved for public release; distribution is unlimited.


\bibliographystyle{unsrt}
\bibliography{report}

\end{document}